\documentclass[a4paper,11pt]{article}

\usepackage[latin1]{inputenc}
\usepackage{amscd,amsmath,amssymb,amsfonts,latexsym,amsthm,mathrsfs}
\usepackage{graphicx,color,framed}
\usepackage{graphicx,color,fancyhdr,framed}%,gensymb,upgreek}
\usepackage{natbib}

\newcommand{\mcg}[2][]{%
\ifthenelse{\equal{#1}{}}{{\mathcal{G}_{#2}}}{\mathcal{G}_{#2}^{(#1)}}%
}
\newcommand{\mcf}[2][]{%
\ifthenelse{\equal{#1}{}}{{\mathcal{F}_{#2}}}{\mathcal{F}_{#2}^{(#1)}}%
}
\newcommand{\mctf}[2][]{%
\ifthenelse{\equal{#1}{}}{\widetilde{\mathcal{F}}_{#2}}{\widetilde{\mathcal{F}}_{#2}^{(#1)}}%
}
\newcommand{\mcbf}[2][]{%
\ifthenelse{\equal{#1}{}}{\overline{\mathcal{F}}_{#2}}{\overline{\mathcal{F}}_{#2}^{(#1)}}%
}

\newcommand{\wrt}{with respect to}
\newcommand{\rhs}{RHS}
\newcommand{\lhs}{LHS}
\newcommand{\iid}{i.i.d.}
\newcommand{\as}{a.s.}

\newcommand{\ie}{i.e.}

\newcommand{\eqsp}{\;}

\newcommand{\rset}{\ensuremath{\mathbb{R}}}
\newcommand{\nset}{\ensuremath{\mathbb{N}}}

\newcommand{\lone}{\ensuremath{\mathsf{L}^1}}
\newcommand{\ltwo}{\ensuremath{\mathsf{L}^2}}
\newcommand{\lp}{\ensuremath{\mathrm{L}^p}}

\newcommand{\1}{\ensuremath{\mathbf{1}}}

\newcommand{\eqdef}{\ensuremath{\stackrel{\mathrm{def}}{=}}}

\newcommand{\rmi}{\ensuremath{\mathrm{i}}}

\newcommand{\esssup}[2][]%
{\ifthenelse{\equal{#1}{}}{\left\| #2 \right\|_\infty}{\left\| #2 \right\|_{#1,\infty}}}
\newcommand{\oscnorm}[2][]%
{\ifthenelse{\equal{#1}{}}{\ensuremath{\operatorname{osc}\left(#2\right)}}{\ensuremath{\operatorname{osc}^{#1}\!\left(#2\right)}}}
\newcommand{\essosc}[3][]%
{\ifthenelse{\equal{#1}{}}{\ensuremath{\operatorname{osc}_{#2}{\left(#3\right)}}}{\ensuremath{\operatorname{osc}^{#1}_{#2}\left(#3\right)}}}

\newenvironment{enum_i}
  {%
  \setlength{\leftmargini}{5ex}%
  \begin{enumerate}}%
  {\end{enumerate}}
  {\end{enumerate}}

\def\LTAR{L}
\def\consistweight{\ensuremath{\mathsf{W}}}
\def\consistfunc{\ensuremath{\mathsf{C}}}
\def\asnormfunc{\ensuremath{\mathsf{A}}}
\def\stdnormfunc{\sigma}
\usepackage{ifthen}
\newboolean{sv_used}
\setboolean{sv_used}{true}
\usepackage{amsthm}
  \newtheorem{thm}{Theorem}
  \newtheorem{prop}[thm]{Proposition}
  \newtheorem{cor}[thm]{Corollary}
  \newtheorem{lem}[thm]{Lemma}
  \newtheorem{defi}{Definition}
  \newtheoremstyle{assumption}
  {6pt}% space above
  {6pt}% space below
  {}%body font
  {-.5em}%indent amount
  {\bfseries}%Theorem head
  {}%punctuation
  {1.5em}%space after theorem head
  {}
  \theoremstyle{assumption}

  \theoremstyle{definition}\newtheorem{rem}{Remark}
  \theoremstyle{definition}

  \newtheoremstyle{example}
  {6pt}% space above
  {6pt}% space below
  {\itshape}%body font
  {-.5em}%indent amount
  {\itshape}%Theorem head
  {}%punctuation
  {1.5em}%space after theorem head
  {}
  \theoremstyle{definition}\newtheorem{ex}{Example}
\newcommand{\PP}{\ensuremath{\operatorname{P}}}
\newcommand{\PE}{\ensuremath{\operatorname{E}}}
\newcommand{\CPE}[3][]
{\ifthenelse{\equal{#1}{}}{\operatorname{E}\left[\left. #2 \, \right| #3 \right]}{\operatorname{E}_{#1}\left[\left. #2 \, \right | #3 \right]}}
\newcommand{\CPP}[3][]
{\ifthenelse{\equal{#1}{}}{\operatorname{P}\left(\left. #2 \, \right| #3 \right)}{\operatorname{P}_{#1}\left(\left. #2 \, \right | #3 \right)}}
\newcommand{\PVar}{\ensuremath{\operatorname{Var}}}

\newcommand{\gauss}{\ensuremath{\operatorname{N}}}

\renewcommand{\mid}{\,|\,}
\newcommand{\ci}[4][]%
{%
\ifthenelse{\equal{#1}{}}{\ensuremath{#2 \perp\!\!\!\perp #3 \mid #4 }}{\ensuremath{#2 \perp\!\!\!\perp #3 \mid #4 \; \: [#1]}}%
}
\newcommand{\dlim}{\ensuremath{\stackrel{\mathcal{D}}{\longrightarrow}}}
\newcommand{\plim}{\ensuremath{\stackrel{\mathrm{P}}{\longrightarrow}}}

\newcommand{\Xset}{\ensuremath{\mathsf{X}}}
\newcommand{\epartset}{\ensuremath{\boldsymbol{\Xi}}}
\newcommand{\epartsigma}{\ensuremath{\mathcal{B}(\boldsymbol{\Xi})}}
\newcommand{\etpartset}{\ensuremath{\boldsymbol{\tilde{\Xi}}}}
\newcommand{\etpartsigma}{\ensuremath{\mathcal{B}(\boldsymbol{\tilde{\Xi}})}}
\newcommand{\Xsigma}[1][]%
{%
\ifthenelse{\equal{#1}{}}{\ensuremath{\mathcal{B}(\Xset)}}{\ensuremath{\mathcal{B}(\Xset^{#1})}}
}

\newcommand{\Yset}{\ensuremath{\mathsf{Y}}}
\newcommand{\Ysigma}{\ensuremath{\mathcal{B}(\Yset)}}

\newcommand{\chunk}[4][]%
{\ifthenelse{\equal{#1}{}}{\ensuremath{{#2}_{#3:#4}}}{\ensuremath{#2^#1}_{#3:#4}}
}
\newcommand{\Q}{\ensuremath{Q}}

\newcommand{\Xinit}{\ensuremath{\chi}}

\newcommand{\URoot}{\ensuremath{R}}
\newcommand{\UCov}[1][]%
{%
\ifthenelse{\equal{#1}{}}{\URoot \URoot^t}{\URoot_{#1} \URoot^t_{#1}}%
}

\newcommand{\VRoot}{\ensuremath{S}}
\newcommand{\VCov}[1][]%
{%
\ifthenelse{\equal{#1}{}}{\VRoot \VRoot^t}{\VRoot_{#1} \VRoot^t_{#1}}%
}

\newcommand{\LDX}[2]{\ensuremath{L}}

\newcommand{\filt}[2][]%
{%
\ifthenelse{\equal{#1}{}}{\ensuremath{\phi_{#2}}}{\ensuremath{\phi_{#1,#2}}}%
}
\newcommand{\pred}[3][]%
{%
\ifthenelse{\equal{#1}{}}{\ensuremath{\phi_{#2|#3}}}{\ensuremath{\phi_{#1,#2|#3}}}%
}
\newcommand{\post}[2][]%
{%
\ifthenelse{\equal{#1}{}}{\ensuremath{\phi_{#2}}}{\ensuremath{\phi_{#1,#2}}}%
}
\newcommand{\logl}[2][]%
{%
\ifthenelse{\equal{#1}{}}{\ensuremath{\ell_{#2}}}{\ensuremath{\ell_{#1,#2}}}%
}
\newcommand{\lhood}[2][]%
{%
\ifthenelse{\equal{#1}{}}{\ensuremath{\mathrm{L}_{#2}}}{\ensuremath{\mathrm{L}_{#1,#2}}}%
}
\newcommand{\cc}[2][]%
{%
\ifthenelse{\equal{#1}{}}{\ensuremath{c_{#2}}}{\ensuremath{c_{#1,#2}}}%
}
\newcommand{\forvar}[2][]%
{%
\ifthenelse{\equal{#1}{}}{\ensuremath{\alpha_{#2}}}{\ensuremath{\alpha_{#1,#2}}}%
}
\newcommand{\nforvar}[2][]%
{%
\ifthenelse{\equal{#1}{}}{\ensuremath{\bar{\alpha}_{#2}}}{\ensuremath{\bar{\alpha}_{#1,#2}}}%
}

\newcommand{\BK}[2][]%
{%
\ifthenelse{\equal{#1}{}}{\ensuremath{\mathrm{\mathrm{B}}_{#2}}}{\ensuremath{\mathrm{B}_{#1,#2}}}%
}
\newcommand{\filtfunc}[2][]%
{%
\ifthenelse{\equal{#1}{}}{\ensuremath{\tau_{#2}}}{\ensuremath{\tau_{#1,#2}}}%
}
\ifthenelse{\boolean{sv_used}}{%
  }
{}

\newcommand{\NISE}[4][]%
{%
\ifthenelse{\equal{#1}{}}{\ensuremath{\tilde{#2}^{\scriptstyle \mathrm{IS}}_{#4}\left(#3 \right)}}{\ensuremath{\tilde{#2}^{\scriptstyle \mathrm{IS}}_{#1,#4}\left(#3 \right)}}%
}
\newcommand{\ISE}[4][]%
{%
\ifthenelse{\equal{#1}{}}{\ensuremath{\widehat{#2}^{\scriptstyle  \mathrm{IS}}_{#4} \left(#3 \right)}}{\ensuremath{\widehat{#2}^{\scriptstyle  \mathrm{IS}}_{#1,#4} \left(#3 \right)}}%
}
\newcommand{\SIRE}[4][]%
{%
\ifthenelse{\equal{#1}{}}{\ensuremath{\hat{#2}^{\scriptstyle  \mathrm{SIR}}_{#4} \left(#3 \right)}}{\ensuremath{\hat{#2}^{\scriptstyle  \mathrm{SIR}}_{#1,#4} \left(#3 \right)}}%
}
\newcommand{\MCE}[3]%
{
{\ensuremath{\hat{#1}^{\scriptstyle  \mathrm{MC}}_{#3} \left(#2 \right)}}%
}
\newcommand{\KUN}{\ensuremath{T}}

\newcommand{\KISS}{\ensuremath{R}}
\newcommand{\KPRI}{\ensuremath{Q}}
\newcommand{\KOPT}{\ensuremath{T^\star}}

\newcommand{\CV}[2]{\mathrm{CV}_{#1}^{(#2)}}
\newcommand{\epart}[2][]
{%
\ifthenelse{\equal{#1}{}}{\ensuremath{\xi_{#2}}}{\ensuremath{\xi_{#2}^{({#1})}}}
}
\newcommand{\etpart}[2][]
{%
\ifthenelse{\equal{#1}{}}{\ensuremath{\tilde{\xi}_{#2}}}{\ensuremath{\tilde{\xi}_{#2}^{({#1})}}}
}
\newcommand{\ebpart}[2][]
{%
\ifthenelse{\equal{#1}{}}{\ensuremath{\bar{\xi}_{#2}}}{\ensuremath{\bar{\xi}_{#2}^{({#1})}}}
}

\newcommand{\etwght}[2][]{%
\ifthenelse{\equal{#1}{}}{\ensuremath{\tilde{\omega}_{#2}}}{\ensuremath{\tilde{\omega}_{#2}^{(#1)}}}}
\newcommand{\ewght}[2][]{%
\ifthenelse{\equal{#1}{}}{\ensuremath{\omega_{#2}}}{\ensuremath{\omega_{#2}^{(#1)}}}}
\newcommand{\funcweight}{\ensuremath{\Phi}}
\newcommand{\sumweight}[2][]{%
\ifthenelse{\equal{#1}{}}{\ensuremath{\Omega_{#2}}}{\ensuremath{\Omega_{#2}^{(#1)}}}}
\newcommand{\tsumweight}[2][]{%
\ifthenelse{\equal{#1}{}}{\ensuremath{\widetilde{\Omega}_{#2}}}{\ensuremath{\widetilde{\Omega}_{#2}^{(#1)}}}}
\newcommand{\indexresidual}[2][]{%
\ifthenelse{\equal{#1}{}}{\ensuremath{J_{#2}}}{\ensuremath{J_{#2}^{(#1)}}}}

\newcommand{\efilt}[2][]%
{%
\ifthenelse{\equal{#1}{}}{\ensuremath{\hat{\phi}_{#2}}}{\ensuremath{\hat{\phi}_{#1,#2}}}%
}
\newcommand{\epost}[3][]%
{%
\ifthenelse{\equal{#1}{}}{\ensuremath{\hat{\phi}_{#2|#3}}}{\ensuremath{\hat{\phi}_{#1,#2|#3}}}%
}

\begin{document}

\begin{center}
{\bf \Large Limit theorems for  weighted samples with applications to Sequential Monte Carlo Methods}\\
\ \\
By R. Douc\footnote{CMAP, \'Ecole Polytechnique, Palaiseau, France. {\sf douc@cmapx.polytechnique.fr}} \hspace{0.05cm} and E. Moulines\footnote{TSI, \'Ecole Nationale Sup\'erieure des T\'el\'ecommunications, Paris, France. {\sf moulines@tsi.enst.fr}}.
\end{center}
\renewcommand{\thefootnote}{}


\begin{abstract}
In the last decade, sequential Monte-Carlo methods (SMC) emer\-ged as
a key tool in computational statistics (see for instance
\cite{doucet:defreitas:gordon:2001}, \cite{liu:2001}, \cite{kuensch:2001}). These algorithms approximate a sequence of
distributions by a sequence of weighted empirical measures
associated to a weighted population of particles. These particles
and weights  are generated recursively according to elementary
transformations: mutation and selection. Examples of
applications include the sequential Monte-Carlo techniques to solve
optimal non-linear filtering problems in state-space models,
molecular simulation, genetic optimization, etc.

Despite many theoretical advances (see for instance
\cite{gilks:berzuini:2001b}, \cite{kuensch:2003}, \cite{delmoral:2004}, \cite{chopin:2004}),
the asymptotic property of these approximations remains of course a question of central
interest. In this paper, we analyze sequential Monte Carlo methods
from an asymptotic perspective, that is, we establish law of large
numbers and invariance principle as the number of particles gets
large. We introduce the concepts of \emph{weighted sample}
consistency and asymptotic normality, and derive conditions under
which the mutation and the selection procedure used in the sequential Monte-Carlo build-up preserve these
properties. To illustrate our findings, we  analyze SMC algorithms
to approximate the filtering distribution in state-space models. We
show how our techniques allow to relax  restrictive technical
conditions used in previously reported works and provide grounds to
analyze more sophisticated sequential sampling strategies.
\end{abstract}
\vspace{0.25cm}\noindent{\bf Short title:} Limit theorems for SMC.

\section{Introduction}
Sequential Monte Carlo (SMC) refer to a class of methods  designed to approximate a \emph{sequence
of probability distributions} over a \emph{sequence of probability space} by a set of points, termed
\emph{particles}  that each have an assigned non-negative weight
and are updated recursively in time. SMC methods can be seen as a combination of the sequential importance sampling introduced method in
\cite{handschin:mayne:1969} and the sampling importance resampling algorithm proposed in \cite{rubin:1987}; it
uses a combination of \emph{mutation} and \emph{selection} steps.
In the \emph{mutation} step, the particles are propagated forward in time using proposal kernels and their importance weights are updated taking into account
the targeted distribution. In the \emph{selection} (or \emph{resampling}) step, particles multiply or die depending on their \emph{fitness} measured
by their importance weights. Many algorithms have been proposed since, which differ in the way the particles and the importance
weights evolve and adapt.

SMC methods have a long history in molecular simulations, where they have been found to be a one of the most
powerful means for the simulation and optimization of chain polymers (see for instance \cite{landau:binder:2000}).
SMC methods have more recently emerged as a key tool  to solve \emph{on-line} prediction / filtering / smoothing problems in a dynamic system. Simple yet flexible SMC methods have been shown to
overcome the numerical difficulties and pitfalls typically encountered with traditional methods based on approximate
non-linear filtering (such as the extended  Kalman filter or gaussian-sum filters);
see for instance \cite{liu:chen:1998}, \cite{liu:2001}, \cite{doucet:defreitas:gordon:2001} and \cite{ristic:arulampalam:gordon:2004} and the references therein.
More recently, SMC methods have been shown to be a promising alternative to Markov Chain Monte Carlo techniques
for sampling complex distributions over large dimensional spaces;
see for instance \cite{gilks:berzuini:2001b} and \cite{Cappe:Guillin:Marin:Robert:2003}.

In this paper, we focus on the asymptotic behavior of the \emph{weighted particle} approximation as
the number of particles tend to infinity. Because the particles interact during the selection steps, they are not independent which make the  analysis of particle approximation a
challenging area of research.  This topic has attracted in recent years a great deal of efforts in recent years making it
a daunting task to give credit every contribution.  The first rigorous convergence result was obtained in \cite{delmoral:1996}, who established the almost-sure convergence
of an elementary SMC algorithm (the so-called \emph{bootstrap filter}). A central limit theorem  for this algorithm was  derived in \cite{delmoral:guionnet:1999}
and refined in \cite{delmoral:miclo:2000}. The proof of the CLT was later simplified  and extended to more general SMC algorithms by
\cite{kuensch:2003} and \cite{chopin:2004}. Bounds on the fluctuations of the particle approximations for different norms were reported
in \cite{crisan:lyons:1997}, \cite{delmoral:miclo:2000} and \cite{crisan:doucet:2002}. \cite{delmoral:2004} provides an up-to-date and thorough
coverage of recent theoretical developments in this area.

With few exceptions (see \cite{chopin:2004} and to a lesser extent, \cite{crisan:doucet:2002} and \cite{kuensch:2003}),
these results apply under simplifying assumptions on the way particles are mutated and selected  which restrict the scope of applicability of
the results only to the most elementary SMC implementations. In particular, all these results assume that selection is performed at
each iteration which implies that the weights are not propagated.
This is clearly an annoying limitations since it has been noticed in practice that resampling the particle system at each time step
is most often not a clever choice. As discussed for example in
\cite[section 2]{liu:chen:1998}, when the weights are nearly constants, resampling  only reduces the number of distinct particles and introduces
extra Monte Carlo variations. Resampling should only be applied when the weights are very skewed:  carrying many particles with
very small importance weights is indeed a waste. Resampling provide chances to \emph{good} particles to replicate and hence rejuvenate the sampler to produce
better \emph{future} particles.

The main purpose of this paper is to derive an asymptotic theory of  weighted system of particles. To our best knowledge, limit theorems
for such weighted  approximations were only considered in \cite{liu:chen:1998}, who mostly sketched consistency proofs. In this paper, we
establish both law of large numbers and central limit theorems, under assumptions that are presumably closed from being minimal. These results
apply not only to the many different implementations of the SMC algorithms, including rather sophisticated schemes such as the
resample and move algorithm \cite{berzuini:gilks:2001} or the auxiliary particle filter by \cite{pitt:shephard:1999}. They cover resampling schedules (when to resample) that can be either
deterministic or dynamic, \ie\ based on the distribution of the importance weights at the current iteration. They also cover sampling scheme that can be either
simple random sampling (with weights) but also residual sampling \citep{liu:chen:1998} or auxiliary sampling \citep{pitt:shephard:1999}. We do not impose
a specific structure on the sequence of the target probability measure; therefore our results apply not only to sequential filtering or smoothing of
state-space contexts, but also to recent algorithms developed for population Monte Carlo or for molecular simulation.

The paper is organized as follows. In section \ref{sec:Notations and Definitions}, we introduce the  definitions of weighted
sample consistency and asymptotic normality; we then discuss the
meaning of these definitions in a simple situation. In section \ref{sec:mainresults}, we  present and discuss the conditions upon
which consistency or / and asymptotic normality of a weighted sample is preserved during the mutation and selection steps.
In section \ref{sec:application to state-space models}, we apply the result to the estimation of the joint smoothing distribution for a
state-space model. In particular, we establish a central limit theorem for a SMC method involving a dynamic resampling scheme. These
results are based on general results on triangular array of martingale increments (in the sense of \cite[Section II.7]{shiryaev:1996}) which
are established in section \ref{sec:WeakLimitsTriangularArray}.

\section{Notations and Definitions}
\label{sec:Notations and Definitions}
All the random variables are defined on a common probability space $(\Omega,\mcf{},\PP)$.
A state space $\Xset$ is said to be \emph{general} if it is equipped with  a countably generated  $\sigma$-field $\Xsigma$.
For a general state space $\Xset$, we denote by $\mathcal{P}(\Xset)$ the set of probability measures on $(\Xset, \Xsigma)$ and
$\mathbb{B}(\Xset)$ (resp. $\mathbb{B}^+(\Xset)$) the set of all $\mathcal{B}(\Xset) / \mathcal{B}(\rset)$-measurable (resp. non-negative) functions from
$\Xset$ to $\rset$ equipped with the Borel $\sigma$-field $\mathcal{B}(\rset)$.
For any   $\mu \in \mathcal{P}(\Xset)$ and $f \in \mathbb{B}(\Xset)$ satisfying $\int_\Xset \mu(dx) |f(x)| < \infty$,
 $\mu(f)$ denotes $\int_\Xset f(x) \mu(dx)$.
Let $\Xset$ and $\Yset$ be two general state spaces. A kernel $V$ from $(\Xset,\Xsigma)$ to $(\Yset,\Ysigma)$ is a
map from $\Xset \times \Ysigma$ into $[0, 1]$ such that for each $A \in \Ysigma$, $x \mapsto V (x,A)$ is a nonnegative
bounded measurable function on $\Xset$ and for each $x \in \Xset$, $A \mapsto V (x, A)$ is a measure
on $\Ysigma$. We say that $V$ is finite if $V(x,\Yset) < \infty$ for any $x \in \Xset$; it is Markovian if $V(x,\Xset) \equiv 1$ for
any $x \in \Xset$.
For any function $ f \in \mathbb{B}(\Xset \times \Yset)$ such that $\int_{\Yset} V(x,dy) |f(x,y)| < \infty$ we denote
by $V(\cdot,f)$  the function
\begin{equation}
\label{eq:nonconventional-notation}
V(\cdot,f): x \mapsto V(x,f) \eqdef \int_{\Yset} V(x,dy) f(x,y) \eqsp.
\end{equation}
The function $V(\cdot,f)$ belongs to $\mathbb{B}(\Xset)$. We sometimes use the abridged notation $Vf$ instead of $V(\cdot,f)$.
For $\nu$ a measure on $(\Xset,\Xsigma)$, we denote by $\nu V$ the measure on
$(\Yset , \Ysigma)$ defined for any $A \in \Ysigma$ by $\nu V(A) = \int_\Xset \nu(dx) V(x,A)$.

%As an example, when considering the estimation of the \emph{filtering} distribution ($i=j=k$ in \eqref{eq:filteringdistribution}),
%we simply set $\Xset$ and $\Xsigma=\Xset$. When estimating the \emph{joint smoothing} distribution ($i=0$, $j=k$ in \eqref{eq:filteringdistribution})
%we set $\Xset= \Xset^{k+1}$ and $\Xsigma= \Xsigma^{\otimes (k+1)}$.
Let $\epartset$ be a general state space, $\mu$ be a probability measure on $(\epartset,\epartsigma)$, $\{ M_N \}_{N \geq 0}$ be a sequence
of integers, and $\consistfunc$ be a subset of $\lone(\epartset,\mu)$.
We approximate the probability measure $\mu$ by points $\epart{N,i} \in \epartset$, $i=1, \dots, M_N$ associated to non-negative weights
$\ewght{N,i} \geq 0$.
\begin{defi}
\label{defi:weighted-sample:consistency}
A weighted sample $\{ (\epart{N,i},\ewght{N,i}) \}_{1 \leq i \leq M_{N}}$ on $\epartset$ is said to be \emph{consistent}
for the probability measure $\mu$ and the set $\consistfunc \subseteq \lone(\epartset,\mu)$ if for any $f \in \consistfunc$, as  $N\to\infty$,
\begin{align}
& \label{eq:defi:weighted-sample:consistency-1}
\sumweight{N}^{-1} \sum_{i=1}^{M_N} \ewght{N,i} f \left(\epart{N,i} \right) \plim \mu(f) \eqsp,\\
& \label{eq:defi:weighted-sample:consistency-2}
\sumweight{N}^{-1} \max_{1 \leq i \leq M_N} \ewght{N,i} \, \plim 0 \quad \mbox{where} \quad \sumweight{N}=\sum_{i=1}^{M_N} \ewght{N,i}\eqsp.
\end{align}
\end{defi}
This definition of weighted sample consistency is similar to the notion of \emph{properly weighted sample} introduced in \cite{liu:chen:1998}.
The difference stems from the condition \eqref{eq:defi:weighted-sample:consistency-2} which implies that the contribution of each individual term in the sum vanishes in the limit as $N \to \infty$, a
condition often referred as \emph{smallness} in the literature.

\begin{ex}[Importance Sampling]
To illustrate the meaning of these conditions, consider the importance sampling estimator.
Let $\mu$ (the \emph{target distribution}) and $\nu$ (the \emph{proposal distribution}) be a \emph{known} (perhaps up to a  normalizing constant)
probability distribution on $(\epartset,\epartsigma)$.
Suppose that $\mu$ is absolutely continuous \wrt\ $\nu$ and denote by $W$ the \emph{importance function},
$W= \beta \frac{d \mu}{d \nu}$ for some $\beta > 0$. Then, the weighted sample
$\left\{ \left(\epart{N,i}, W(\epart{N,i}) \right) \right \}_{1 \leq i \leq N}$, where $\{ \epart{N,i} \}_{1 \leq i \leq N}$
are \iid\ $\nu$-distributed is consistent for $\{\mu, \lone(\epartset,\mu) \}$. Eq.~\eqref{eq:defi:weighted-sample:consistency-1} follows
from the law of large numbers. For any $\epsilon, C > 0$,
\[
 N^{-1}  \PE\left[\max_{1 \leq i \leq N} \ewght{N,i} \right] \leq
 \epsilon  + N^{-1} \sum_{i=1}^{N} \PE \left[ \ewght{N,i} \1_{\{ \ewght{N,i} \geq \epsilon N \}} \right]
\leq \epsilon + \nu \left( W \1 {\{ W \geq C \}} \right)
\]
for sufficiently large $N$. Because $\nu(W) = \beta$, $\nu \left( W \1 {\{ W \geq C \}} \right)$ converges to zero as $C \to \infty$,
and \eqref{eq:defi:weighted-sample:consistency-2} follows because $\sumweight{N}/N \plim \beta$.
\end{ex}
Of course, importance sampling is only one technique of the many techniques that can be used to obtain a consistent weighted samples;
other approaches will be considered below. For more complex sampling schemes $\consistfunc$ can be a proper subset of $\lone(\epartset,\mu)$.
In order to obtain sensible results we  restrict our attention to classes of sets which are sufficiently rich.
\begin{defi}[Proper Set]
A set $\consistfunc$ of real-valued measurable functions on $\epartset$ is said to be \emph{proper} if the following conditions are satisfied.
\begin{enum_i}
\item $\consistfunc$ is a linear space: for any $f$ and $g$ in $\consistfunc$ and reals $\alpha$ and $\beta$, $\alpha f + \beta g \in \consistfunc$ \eqsp,
\item If $g \in \consistfunc$ and $f$ is measurable with $|f| \leq |g|$, then $|f| \in \consistfunc$\eqsp,
\item For all $c$, the constant function $f \equiv c$ belongs to $\consistfunc$.
\end{enum_i}
\end{defi}
For any function $f$, define the positive and negative parts of it by
$$
f^+ \eqdef f \vee 0  \quad \text{and} \quad f^- \eqdef (-f) \vee 0 \eqsp,
$$
and note that $f^+$ and $f^-$ are both dominated by $|f|$. Thus,
if $|f| \in \consistfunc$ then
$f^+$ and $f^-$ both belong to $\consistfunc$
and so does $f= f^+ - f^-$.
It is easily seen that for any $p \geq 0$ and any measure $\mu$ on $(\epartset,\epartsigma)$, the set $\lp(\epartset,\mu)$ is proper.

A classical way to strengthen consistency it is to consider distributional convergence of the normalized difference.
Recall that (see \cite{aldous:eagleson:1978}, \cite[chapter 3]{hall:heyde:1981}) if a sequence of rv's converges in distribution to $X$, (written $X_N \dlim X$),
the convergence is said to be \emph{stable} (written $X_N \dlim X$ (stably)), if for any set $B \in \mcf{}$ and for a countable dense set of points $x \in \rset$,
\[
\PP( X_{N} \leq x, B) \quad \text{exists}.
\]
In other word, for all events $B \in \mcf{}$ such that $\PP(B) > 0$, the distribution of $X_N$ conditional on $B$ converges in distribution
to some (proper) distribution which may depend on $B$.
The convergence is \emph{mixing} (written $X_N \dlim X$ (mixing)) if in a stable limit result $X_N \dlim X^*$ (stably) the limit random variable $X^*$
can be chosen to be independent of $\mathcal{F}$. Mixing is in particular useful when studying random sum limit theorems.

Let $\mu$  be a probability measure  on $(\epartset,\epartsigma)$,
$\gamma$ be a finite measure on $(\epartset, \epartsigma)$, $\asnormfunc  \subseteq \lone(\epartset,\mu)$ and
$\consistweight \subseteq \lone(\epartset,\gamma)$ be sets of real-valued measurable functions on $\epartset$,
$\stdnormfunc$ a real non negative function on $\asnormfunc$,  and $\{a_N \}$ be a non-decreasing real sequence diverging to infinity.

\begin{defi}
A weighted sample $\{ (\epart{N,i},\ewght{N,i}) \}_{1 \leq i \leq M_N}$ on $\epartset$ is said to be  \emph{asymptotically normal} for $(\mu,\asnormfunc,\consistweight,\stdnormfunc,\gamma,\{a_N\})$
if
\begin{align}
& a_N \sumweight{N}^{-1} \sum_{i=1}^{M_N} \ewght{N,i} \{ f(\epart{N,i}) - \mu(f) \}
\dlim \gauss \{ 0, \sigma^2(f) \} \quad \text{(mixing)} \quad && \text{for any $f \in \asnormfunc$} \eqsp,\label{eq:AN1}\\
& a_N^2 \sumweight{N}^{-2} \sum_{i=1}^{M_N} \ewght{N,i}^2 f(\epart{N,i}) \plim \gamma(f) \quad && \text{for any $f \in \consistweight$} \label{eq:AN2}\\
&  a_N \sumweight{N}^{-1} \max_{1\leq i \leq M_N}\ewght{N,i} \plim 0\label{eq:AN3}
\end{align}
\end{defi}

We stress that the rate $\{a_N\}$ can be different from $\sqrt{M_N}$ because of the dependence among the weighted sample
$\{ (\epart{N,i}, \ewght{N,i}) \}_{1 \leq i \leq M_N}$ introduced by the different transformations undergone by the weighted sample.

\begin{ex}[Importance Sampling (continued)]
The rationale for the conditions \eqref{eq:AN2} and \eqref{eq:AN3} is best understood by considering again the importance sampling
example. Assume that $\nu(W^2) < \infty$, where as above $W= \beta \; \frac{d \mu}{d \nu}$.
Define the subsets $\asnormfunc \subset \mathbb{B}(\epartset)$ and $\consistweight \subset \mathbb{B}(\epartset)$,
\[
\asnormfunc \eqdef \left\{ f \in \mathbb{B}(\epartset)\, : \nu(W^2 f^2) < \infty \right\} \quad \text{and} \quad
\consistweight \eqdef \left\{ f \in \mathbb{B}(\epartset) \, : \nu(W^2 |f|) < \infty \right\} \eqsp.
\]
For $f \in \asnormfunc$, define $S_{N}(f)  \eqdef \sumweight{N}^{-1} \sum_{i=1}^{N} \ewght{N,i} [f(\epart{N,i})- \mu(f)]$.
Using $\sumweight{N} / N \plim \beta$ and  \cite[Theorem 2]{aldous:eagleson:1978}, $N^{1/2} S_{N}(f) \dlim S$ (mixing),
where $S$ is Gaussian random variable  with zero-mean and variance
\begin{equation}
\label{eq:DefinitionVarianceAS}
\sigma^2(f) \eqdef \nu \left\{ \left( \frac{d \mu}{d \nu} \right)^2 [ f - \mu(f)]^2 \right\} \eqsp.
\end{equation}
The law of large numbers for \iid\ sequence implies that for any $f \in \consistweight$,
\[
N \sumweight{N}^{-2} \sum_{i=1}^{N} \ewght{N,i}^2 f(\epart{N,i}) = N \sumweight{N}^{-2} \sum_{i=1}^{N} W^2(\epart{N,i}) f(\epart{N,i}) \plim  \gamma(f)= \nu \left[\left( \frac{d \mu}{d \nu} \right)^2 f \right]\eqsp,
\]
showing \eqref{eq:AN2}. For any $\epsilon, C > 0$,
\[
 N^{-1}   \PE\left[\max_{1 \leq i \leq N} \ewght{N,i}^2 \right] \leq \epsilon^2 +
 N^{-1} \sum_{i=1}^{N} \PE \left[ \ewght{N,i}^2 \1_{\{ \ewght{N,i}^2 \geq \epsilon^2 N  \}} \right] \\
\leq \epsilon^2 + \nu \left( W^2 \1 {\{ W^2 \geq C \}} \right)
\]
for sufficiently large $N$. Because $\asnormfunc$ is a proper set,  $\nu \left( W^2 \1 { \{ W \geq C \}} \right)$
go to zero as $C \to \infty$. Using again $\sumweight{N} / N \plim \beta$, the previous display implies that
\begin{equation}
\label{eq:UniformSmallnessCondition}
N^{1/2} \sumweight{N}^{-1} \max_{1 \leq i \leq N} \ewght{N,i} \plim 0 \eqsp,
\end{equation}
showing \eqref{eq:AN3}: any individual term in the sum is small compared to the square root of the variance, a
condition which is known to be necessary for a central limit theorem to hold for triangular arrays of independent random variables, see
\cite{petrov:1995}.
%If $\kappa_N$ is a sequence of integer-valued random variables, and if $\kappa_N / N$ converges in probability to a positive random variable
%$M$, then $\kappa_N^{1/2} S_{\kappa_N} \dlim \mathcal{N}(0,\sigma^2(f))$ (see \cite{csorgo:fischler:1971}).
Stability in the limit can be used to deal with situations where the number of terms in the sum are random (see \cite{csorgo:fischler:1971}). This is considered
in a companion paper.
\end{ex}

\section{Main results}
\label{sec:mainresults}
To analyze the  sequential Monte Carlo methods discussed in the introduction, we now need to
study how the  mutation step and  the resampling step affects the consistent or / and
asymptotically normal  weighted sample.
\subsection{Mutation}
\label{subsec:mutation}
To study SISR algorithms, we need first to show that
when moving the particles using a Markov kernel and then assigning them appropriately defined importance weights, we transform a weighted sample consistent  (or asymptotically normal)
for a distribution $\nu$ on a general state space $(\epartset,\epartsigma)$ into a weighted sample consistent (or asymptotically normal) for a distribution $\mu$
on $(\etpartset,\etpartsigma)$. Let $\LTAR$ be a kernel from $(\epartset,\epartsigma)$ to $(\etpartset,\etpartsigma)$ such that
$\nu \LTAR(\etpartset) > 0$ and for any $f \in \mathbb{B}(\etpartset)$,
\begin{equation}
\label{eq:definition-mu}
\mu(f)= \frac{\nu \LTAR (f)}{\nu \LTAR(\etpartset)} = \frac{\iint \nu( d \epart{}) \LTAR(\epart{},d \etpart{}) f( \etpart{})}{\iint \nu( d \epart{}) \LTAR(\epart{},\etpartset)} \eqsp.
\end{equation}
There exist of course many such kernels: one may set for example $\LTAR(\epart{},\tilde{A})= \mu(\tilde{A})$ for any $\epart{} \in \epartset$,
but, as we will see below, this is not usually the most appropriate choice.

We wish to transform a weighted sample $\{(\epart{N,i},\ewght{N,i})\}_{1 \leq i \leq M_N}$ targeting the distribution
$\nu$ on $(\epartset,\epartsigma)$ into a weighted sample $\{ (\etpart{N,i},\etwght{N,i}) \}_{1 \leq i \leq \tilde{M}_N}$ targeting $\mu$
on $(\etpartset,\etpartsigma)$. The number of particles $\tilde{M}_N$ is set to be  $\tilde{M}_N = \alpha_N M_N$ where
$\alpha_N$ is the number of offsprings of each particle.
The use of multiple offspring for a particle has been suggested in the context of SIR by \cite{rubin:1987}: when the mutation step is followed by a
resampling step, an increase in the number of distinct particles in the mutation step
will typically increase the number of distinct particles \textbf{after} the resampling step. In the sequential context,
this operation is regarded as  a practical mean for contending particle impoverishment. These offsprings are proposed using  Markov kernels
$\KISS$, $k=1, \dots, \alpha_N$, from $(\epartset,\epartsigma)$ to $(\etpartset,\etpartsigma)$.
Implicitly, we suppose that  sampling  from the proposal kernels $\KISS$ is doable.

Most importantly, we assume that for any $\epart{} \in \epartset$,
the probability measure $\LTAR(\epart{},\cdot)$ on $(\etpartset,\etpartsigma)$   is absolutely continuous \wrt\  $\KISS$,
which we denote $\LTAR(\epart{},\cdot) \ll \KISS(\epart{},\cdot)$ and define
\begin{equation}
\label{eq:definitionW}
W(\epart{},\etpart{}) \eqdef \frac{d \LTAR(\epart{},\cdot)}{d \KISS(\epart{},\cdot)}(\etpart{})
\end{equation}
%To simplify the discussion, it is assumed in the sequel that
%the number of offsprings $\alpha_N$ is either
%\begin{itemize}
%\item constant, $\alpha_N= \alpha$, for all $N$ where $\alpha$ is an integer; in this case, we define the \emph{average kernel}
%defined by
%\begin{equation}
%\label{eq:average-kernel}
%\KISS \eqdef \alpha^{-1} \sum_{k=1}^{\alpha} \KISS_k \eqsp,
%\end{equation}
%\item $\lim_{N \to \infty} \alpha_N$ exists and is equal to $\alpha$ (which is possibly infinite).
%In the latter case, we suppose that $\KISS_k \equiv \KISS$ for all $k$,
%where $\KISS$ is a Markov kernel from $(\epartset,\epartsigma)$ to $(\etpartset,\etpartsigma)$.
%\end{itemize}
%It is possible to relax the latter assumption, but at the expense of substantial technicalities.
The new weighted sample $\{(\etpart{N,i},\etwght{N,i})\}_{1 \leq i \leq \tilde{M}_N}$ is constructed as follows.
We draw new particle positions $\{ \etpart{N,j} \}_{1 \leq j \leq \tilde{M}_N}$ conditionally independently given
$\mcf{N,0} =  \sigma \left( \{ (\epart{N,i}, \ewght{N,i}) \}_{1 \leq i \leq M_N} \right)$
with distribution given for $i=1, \dots, M_N$, $k=1, \dots, \alpha_N$  and $A \in \etpartsigma$ by
\begin{equation}
\label{eq:mutationstep-1}
\CPP{\etpart{N,\alpha_N (i-1) + k} \in A}{\mcf{N,0}} = \KISS(\epart{N,i},A)
\end{equation}
and associate to each new particle positions  the importance weight:
\begin{equation}
\label{eq:mutationstep-3}
\etwght{N,\alpha_N (i-1) + k} = \ewght{N,i} W(\epart{N,i},\etpart{N,\alpha_N (i-1) + k}) \eqsp, \ \text{where}\  i= 1, \dots, M_N \eqsp, \ k=1, \dots, \alpha_N \eqsp.
\end{equation}
The mutation step is \emph{unbiased} in the sense that, for any $f \in \mathcal{B}(\etpartset)$ and $i=1, \dots, M_N$,
\begin{equation}
\label{eq:unbiased-mutation}
\sum_{j=\alpha_N (i-1)+1}^{\alpha_N i} \CPE{\etwght{N,j} f(\etpart{N,j})}{\mathcal{F}_{N,j-1}}=
\alpha_N \ewght{N,i} \LTAR(\epart{N,i},f) \eqsp,
\end{equation}
where for $j=1, \dots, \tilde{M}_N$, $\mcf{N,j} \eqdef \mcf{N,0} \vee \sigma (\{ \etpart{N,l} \}_{1 \leq l \leq j})$.
The following theorems state conditions under which the mutation step described above preserves the weighted sample consistency.
\begin{thm}
\label{thm:LLN-mutation}
Let $\nu$ be a probability measure on $(\epartset,\epartsigma)$, $\mu$ be a probability measure on $(\etpartset,\etpartsigma)$.
Let $\LTAR$ be a finite kernel from $(\epartset,\epartsigma)$ to $(\etpartset,\etpartsigma)$ satisfying $\nu \LTAR(\etpartset) > 0 $
any $\tilde{A} \in \etpartsigma$, $\mu(\tilde{A})= \nu \LTAR(\tilde{A}) / \nu \LTAR(\etpartset)$.
Assume that
\begin{itemize}
\item[(i)] the weighted sample $\{ (\epart{N,i},\ewght{N,i}) \}_{1 \leq i \leq M_N}$ is consistent for
$(\nu,\consistfunc)$, where  $\consistfunc \subset \lone(\epartset,\nu)$ is a proper set \eqsp,
\item[(ii)]  the function $\LTAR(\cdot,\etpartset)$ belongs to $\consistfunc$ \eqsp,
\item[(iii)] for any $\epart{} \in \epartset$, $ \LTAR(\epart{},\cdot) \ll \KISS(\epart{},\cdot) $ \eqsp.
\end{itemize}
Then,
$\tilde{\consistfunc}$ is a proper set and the weighted sample $\{ (\etpart{N,i},\etwght{N,i}) \}_{1 \leq i \leq \tilde{M}_N}$ defined by
\eqref{eq:mutationstep-1} and \eqref{eq:mutationstep-3}
is consistent for $(\mu,\tilde{\consistfunc})$, where
\begin{equation}
%\label{eq:definition:tildeconsistfunc}
\tilde{\consistfunc} \eqdef \{ f \in \lone(\etpartset,\mu), \, \LTAR(\cdot,|f|) \in \consistfunc \} \eqsp,
\end{equation}
\end{thm}
We now turn to prove the asymptotic normality.
\begin{thm}
\label{thm:CLT-mutation}
Suppose that the assumptions of Theorem \ref{thm:LLN-mutation} hold.
Let $\gamma$ be a finite measure on $(\epartset, \epartsigma)$,
$\asnormfunc \subseteq \lone(\epartset,\nu)$, and $\consistweight \subseteq \lone(\epartset,\gamma)$ be proper sets,
$\stdnormfunc$ be a non negative function on $\asnormfunc$. Assume in addition that
\begin{enum_i}
\item the weighted sample $\{ (\epart{N,i},\ewght{N,i}) \}_{1 \leq i \leq M_N}$ is  asymptotically normal for
$(\nu,$ $\asnormfunc,$ $\consistweight,$ $\stdnormfunc,$ $\gamma,$ $\{M^{1/2}_N\} )$~\eqsp,
\item the function $\KISS(\cdot,W^2)$ belongs to $\consistweight$\eqsp.
\end{enum_i}
Then, if $\alpha_N \to \alpha$,  the weighted sample $\{ (\etpart{N,i},$ $\etwght{N,i}) \}_{1 \leq i \leq \tilde{M}_N}$ is asymptotically normal for $(\mu,$ $\tilde{\asnormfunc},$ $\tilde{\consistweight}$,  $\tilde{\stdnormfunc}$, $\tilde{\gamma},$ $\{M_N^{1/2}\})$ with
\begin{align*}
&\tilde{\asnormfunc}  \eqdef  \left\{ f \in \lone(\etpartset,\mu) \, :  \LTAR (\cdot,|f|) \in \asnormfunc \quad \text{and} \quad
\KISS(\cdot, W^2 f^2)  \in \consistweight  \right \} \eqsp, \\
&\tilde{\consistweight} \eqdef \left\{ f \in \lone(\etpartset,\mu) \, :   \KISS(\cdot, W^2 |f|) \in \consistweight \right\} \eqsp, \\
& \tilde{\stdnormfunc}^2(f) \eqdef \stdnormfunc^2\{  \tilde{\LTAR}[f - \mu(f)] \} + \alpha^{-1} \gamma \KISS\{[\tilde{W}f - \KISS(\cdot,\tilde{W}f)]^2\} \eqsp,
&& \text{for all $f \in \tilde{\asnormfunc}$}\eqsp,\\
& \tilde{\gamma}(f)\eqdef \alpha^{-1}\gamma\KISS(\tilde{W}^2 f) \eqsp,
&& \text{for all  $f \in \tilde{\consistweight}$}\eqsp,
\end{align*}
where $\tilde{\LTAR} \eqdef [\nu \LTAR(\etpartset)]^{-1} \LTAR$ and $\tilde{W}= [\nu \LTAR(\etpartset)]^{-1} W$. Moreover, $\tilde{\asnormfunc}$ and $\tilde{\consistweight}$ are proper sets.
\end{thm}
\begin{rem}
The use of different proposal kernels $\KISS_k$ is sometimes advisable, to add variety in the type of moves that are applied to the particle system.
The idea of using multiple type of moves is exploited in particular in adaptive importance
sampling (see \cite{douc:guillin:marin:robert:2004} for an illustration). The results above can be adapted directly when the number of
offsprings $\alpha_N= \alpha$ does not depend on $N$. To handle this case, define the average kernel $\KISS = \alpha^{-1} \sum_{k=1}^\alpha \KISS_k$
and suppose that for any $\epart{} \in \epartset$, $\LTAR(\epart{},\cdot) \ll \KISS(\epart{},\cdot)$ (note that the individual proposal kernels
need not be absolutely continuous \wrt\ $\LTAR$).
For $i=1, \dots, M_N$, $k=1, \dots, \alpha_N$  and $A \in \etpartsigma$, we draw the new particle positions according to
\begin{equation}
\label{eq:mutationstep-1-multiple}
\CPP{\etpart{N,\alpha_N (i-1) + k} \in A}{\mcf{N,0}} = \KISS_k(\epart{N,i},A) \eqsp .
\end{equation}
We associate to each particle the importance weight,
\begin{equation}
\label{eq:mutationstep-3-multiple}
\etwght{N,\alpha_N (i-1) + k} = \ewght{N,i} W(\epart{N,i},\etpart{N,\alpha_N i + k}) \eqsp, \quad \text{where}\  i= 1, \dots, M_N \eqsp, \ k=1, \dots, \alpha_N \eqsp.
\end{equation}
where $W$ is defined as in \eqref{eq:definitionW}. With these definitions and notations, the results of Theorems \ref{thm:LLN-mutation} and \ref{thm:CLT-mutation} continue to
hold.
\end{rem}

To illustrate the meaning of these statements, consider again the elementary example of importance sampling. It turns out that reweighting the particles without moving them is a particular case of mutation. Thus, Theorem \ref{thm:LLN-mutation} and \ref{thm:CLT-mutation} may apply in this context.
\begin{ex}
\label{ex:importance-sampling}
Let $\nu$ and $\mu$ be two probability measures on $(\epartset,\epartsigma)$. Assume that $\mu \ll \nu$. Set $\LTAR(\epart{},\cdot)= \beta \frac{d \mu}{d \nu}(\epart{}) \delta_{\epart{}}$,
for some $\beta > 0$, where $\delta_{\epart{}}$ is the Dirac measure. With this definition of $\LTAR$, we obviously have $\mu(\cdot)= \nu \LTAR(\cdot)/ \nu \LTAR(\epartset)$.
If $\{ \epart{N,i} \}_{1 \leq i \leq M_N}$ are independent $\nu$-distributed random variables, then the weighted sample $\{ (\epart{N,i},1) \}_{1 \leq i \leq M_N}$ is consistent for
$\left(\nu, \lone(\epartset,\nu) \right)$. Put for $\epart{} \in \epartset$ and $A \in \epartsigma$, $\KISS(\epart{},A) = \delta_{\epart{}}(A)$. Then,
Theorem \ref{thm:LLN-mutation} shows that $\{ (\etpart{N,i}, W(\epart{N,i},\etpart{N,i})) \}_{1 \leq i \leq M_N}$
where $W(\epart{},\etpart{}) = \beta \frac{d \mu}{d \nu}(\etpart{})$
is consistent for $\{\mu, \lone(\epartset,\mu)\}$. We now turn to check the asymptotic normality. From the CLT for \iid\ sequence,
the weighted sample is $\{ (\epart{N,i}, 1)\}_{1 \leq i \leq M_N}$ is asymptotically normal for $(\nu,$ $\ltwo(\Xset,\nu),$ $\lone(\Xset,\nu)$,
$\sigma^2,$ $\nu,$ $M_N^{1/2}\}$, where for $f \in \asnormfunc \eqdef \ltwo(\Xset,\nu)$, $\nu(f) = \nu\{ [f - \nu(f)]^2 \}$.
Suppose that $\nu(W^2) < \infty$. Because $\KISS(\cdot,W^2)= W^2 \in \consistweight = \lone(\Xset,\nu)$, then
Theorem \ref{thm:CLT-mutation} shows that $\{ (\etpart{N,i}, W(\epart{N,i}, \etpart{N,i})) \}_{1 \leq i \leq M_N}$ is asymptotically normal
for $(\mu,$ $\tilde{\asnormfunc},$ $\tilde{\consistfunc},$ $\tilde{\stdnormfunc},$ $\tilde{\gamma},$ $\sqrt{M_N} \}$,
where $\tilde{\asnormfunc}= \{ f \in \lone(\epartset,\mu), \nu(W^2f^2) < \infty \}$,
$\tilde{\consistweight} = \{ f \in \lone(\epartset,\mu), \nu(W^2|f| ) < \infty \}$, $\tilde{\gamma}(f) =  \nu \left[ W^2 f \right]$ for
$f \in \tilde{\consistweight}$, and $\tilde{\sigma}^2(f)= \PVar_\nu \left[ W f \right]$ for  $f \in \tilde{\asnormfunc}$. These
expressions are particularly simple because, in this case, $Wf = \KISS(\cdot,Wf)$, which is of course false in general.
\end{ex}

\subsection{Resampling}
\label{subsec:resampling}
Resampling is the second basic type of transformation used in sequential Monte-Carlo methods;
resampling converts a weighted sample $\{ (\epart{N,i},\ewght{N,i}) \}_{1 \leq i \leq M_N}$
targeting a distribution $\nu$ into an equally weighted sample $\{ (\etpart{N,i},1)\}_{1 \leq i \leq \tilde{M}_N}$
targeting the same distribution $\nu$ on $(\epartset, \epartsigma)$. Note that the number of resampled particles, denoted $\tilde{M}_N$, might well
differ from the initial number of particles $M_N$. In this section, we  will focus on importance sampling estimators
that satisfy the following \emph{unbiasedness condition}:
\begin{equation}
\label{eq:unbiased-sampling-condition}
\text{for any $f \in \mathbb{B}(\epartset)$}, \quad
\CPE{\tilde{M}_N^{-1} \sum_{i=1}^{\tilde{M}_N} f(\etpart{N,i})}{\mcf{N}} = \sumweight{N}^{-1} \sum_{i=1}^{M_N} \ewght{N,i} f(\epart{N,i}) \eqsp,
\end{equation}
where $ \mcf{N}= \sigma(\{ (\epart{N,i},\ewght{N,i}) \}_{1 \leq i \leq M_N} $. The unbiasedness condition \eqref{eq:unbiased-sampling-condition} implies that resampling can only increase the variability of the importance sampling estimator
and has thus, in one step,  an adversary effect.
In the sequential context however, resampling is nevertheless essential by  removing
particles with small importance weights and producing multiple copies of particles with large importance weights,
that help in generating better future samples. There are many different unbiased resampling procedures
described in the literature. The most elementary procedure is the so-called \emph{multinomial} sampling in which we
draw conditionally independently given $\mcf{N}$ integer-valued random variables $\{ I_{N,k} \}_{1 \leq k \leq \tilde{M}_N}$ with distribution
\begin{equation}
\label{eq:multinomial-sampling-1}
\CPP{I_{N,k}=i}{\mcf{N}}= \sumweight{N}^{-1} \ewght{N,i} \eqsp, \quad i = 1, \dots, M_N \eqsp,
\end{equation}
and set
\begin{equation}
\label{eq:multlinomial-sampling-2}
\etpart{N,i}= \epart{N,I_{N,i}} \eqsp, \quad \text{for} \  i=1, \dots, \tilde{M}_N \eqsp.
\end{equation}
To bring down the extra-variability incurred by multinomial sampling, other resampling strategies have been considered in the literature. A possible solution, investigated in \cite{liu:chen:1995},
consists in  using a combination of deterministic plus residual sampling. This scheme consists in retaining at
least $\lfloor \sumweight{N}^{-1} \tilde{M}_N \ewght{N,i}  \rfloor$, $i=1, \dots, M_N$ copies of the
particles and then reallocating the remaining particles by applying the multinomial sampling procedure but on the residual
importance weights, $\tilde{M}_N \sumweight{N}^{-1}  \ewght{N,i} - \lfloor \tilde{M}_N \sumweight{N}^{-1}  \ewght{N,i}  \rfloor$.
Define $\indexresidual{N,0}=0$ and for $k \geq 1$, $\indexresidual{N,k} = \sum_{i=1}^{k} \lfloor  \tilde{M}_N \sumweight{N}^{-1} \ewght{N,i} \rfloor$. The residual resampling algorithm is defined as follows:
\begin{itemize}
\item [(i)] For $i=1, \dots, M_N$, assign
\begin{equation}
\label{eq:multinomial-sampling-resid-1}
\etpart{N,j}=\epart{N,i} \eqsp, \quad \text{for all $  j= \indexresidual{N,i-1}+1, \dots,  \indexresidual{N,i}$}\eqsp.
\end{equation}
\item [(ii)] Draw, conditionally independently given
$ \mcf{N}= \sigma(\{ (\epart{N,i},\ewght{N,i}) \}_{1 \leq i \leq M_N} $
random variables $\{ I_{N,k} \}_{1 \leq k\leq \tilde{M}_N - \bar{M}_N}$ where
\begin{align}
&\bar{M}_N \eqdef \indexresidual{N,M_N} = \sum_{i=1}^{M_N} \lfloor  \tilde{M}_N \sumweight{N}^{-1} \ewght{N,i} \rfloor\eqsp,  \nonumber \\
&\CPP{I_{N,k}=i}{\mcf{N}}= \frac{\tilde{M}_N \sumweight{N}^{-1} \ewght{N,i}- \lfloor \tilde{M}_N \sumweight{N}^{-1} \ewght{N,i}  \rfloor}{\tilde{M}_N-\bar{M}_N} \eqdef \etwght{N,i} \eqsp,\quad  i = 1, \dots, M_N \eqsp.\label{eq:multinomial-sampling-resid-2}
\end{align}
and set $\etpart{N,\bar{M}_N+k}=\epart{N,I_{N,k}}$ for all $1 \leq k \leq \tilde{M}_N - \bar{M}_N$.
\end{itemize}

If the weighted sample $\{ (\epart{N,i},\ewght{N,i}) \}_{1 \leq i \leq M_N}$ is consistent for $(\nu,\consistfunc)$, where
$\consistfunc$ is a proper subset of $\mathbb{B}(\Xset)$,
it is a natural question to ask whether $\{ (\etpart{N,i}, 1) \}_{1 \leq i \leq \tilde{M}_N}$ is consistent for $\nu$
and, if so, what an appropriately defined class of functions on $\epartset$ might be. It happens that a fairly general result can be obtained
in this case.
\begin{thm}
\label{thm:LLN-unbiased-sampling}
Let $\nu$ a probability measure on $(\epartset, \epartsigma)$ and
$\consistfunc \subseteq \lone(\epartset,\nu)$ be a proper set of functions.
Assume that the weighted sample $\{ (\epart{N,i},$ $\ewght{N,i}) \}_{1 \leq i \leq M_N}$ is consistent for $(\nu,\consistfunc)$,
where $\consistfunc \subseteq \mathbb{B}(\epartset)$ is a proper set.
Then, the uniformly weighted sample $\{ (\etpart{N,i},$ $ 1)  \}_{1 \leq i \leq \tilde{M}_N}$ obtained
using \textbf{any} unbiased resampling scheme (\ie\ satisfying \eqref{eq:unbiased-sampling-condition})  is consistent for $(\nu,\consistfunc)$.
\end{thm}
It is also sensible in this discussion to strengthen the requirement of consistency into asymptotic normality, and
again prove that the sampling operation transform an asymptotically normal weighted sample for $\nu$ into an asymptotically normal sample for $\nu$
(for appropriately defined class of functions, normalizing factors, etc.). We consider first the multinomial sampling algorithm.
\begin{thm}
\label{thm:CLT-independent-sampling}
Suppose that the assumptions of Theorem \ref{thm:LLN-unbiased-sampling} hold.
Let $\gamma$ be a finite measure on $(\epartset, \epartsigma)$,
$\asnormfunc \subseteq \lone(\epartset,\nu)$ and $\consistweight \subseteq \lone(\epartset,\gamma)$ be proper sets,
$\stdnormfunc$ be a non negative function on $\asnormfunc$ and $\{a_N\}$ be a non-negative sequence. Define
\begin{equation} \label{eq:definition-asnormfunc-sampling}
\tilde{\asnormfunc} \eqdef \left\{ f \in \asnormfunc, \,  f^2  \in \consistfunc \right \} \eqsp,
\end{equation}
Assume in addition that
\begin{enum_i}
\item the weighted sample $\{ (\epart{N,i},$ $ \ewght{N,i}) \}_{1 \leq i \leq M_N}$ is
asymptotically nor\-mal for $(\nu,$ $\asnormfunc,$ $\consistweight$, $\stdnormfunc,$ $\gamma,$ $\{a_N\})$
\item $ \lim_{N \to \infty} a_N = \infty$ and $\lim_{N \to \infty} a^2_N / \tilde{M}_N= \beta$, where $\beta \in [0, \infty]$.
\end{enum_i}
Then $\tilde{\asnormfunc}$ is a proper set and the following holds true for the resampled system
$\{(\etpart{N,i},1)\}_{1 \leq i\leq \tilde{M}_N}$ defined as in \eqref{eq:multinomial-sampling-1} and \eqref{eq:multlinomial-sampling-2}.
\begin{enum_i}
\item If $\beta < 1$, then $\{(\etpart{N,i},1)\}_{1 \leq i \leq \tilde{M}_N}$ is
asymptotically normal for $(\nu,$ $\tilde{\asnormfunc},$ $\consistfunc,$ $\tilde{\stdnormfunc},$ $\tilde{\gamma},$ $\{ a_N \})$ with
\begin{align*}
& \tilde{\stdnormfunc}^2(f)=\beta \PVar_\nu(f)+ \stdnormfunc^2(f) \quad \text{for any $f \in \tilde{\asnormfunc}$}\eqsp,\\
  &\tilde{\gamma}= \beta \nu \eqsp.
\end{align*}
\item If $\beta \geq 1$, then $\{(\etpart{N,i},1)\}_{1 \leq i \leq \tilde{M}_N}$ is
 asymptotically normal for $(\nu,$ $\tilde{\asnormfunc},$ $\consistfunc,$ $ \tilde{\stdnormfunc},$ $\tilde{\gamma},$ $\{ \tilde{M}^{1/2}_N \})$ with
 \begin{align*}
& \tilde{\stdnormfunc}^2(f)= \PVar_\nu(f)+  \beta^{-1}\stdnormfunc^2(f) \quad \text{for any $f \in \tilde{\asnormfunc}$}\eqsp,\\
   &\tilde{\gamma}= \nu\eqsp.
 \end{align*}
\end{enum_i}
\end{thm}
\begin{ex}
Let $\nu, \mu \in \mathcal{P}(\epartset)$ and suppose that $\mu \ll \nu$.
As shown in the preceding example, if $\{ \epart{N,i} \}_{1 \leq i \leq M_N}$ is an \iid\ sample from  $\nu$,
then the weighted sample $\left\{ (\epart{N,i}, W(\epart{N,i})) \right\}_{1 \leq i \leq M_N}$, where $W \propto d \mu / d \nu$
is consistent for $\left\{ \mu,\lone(\epartset,\mu) \right\}$ and asymptotically normal for $(\mu,$ $\asnormfunc,$ $\consistweight,$ $\stdnormfunc,$ $\gamma,$ $\sqrt{M_N} \}$,
where $\asnormfunc= \{ f \in \lone(\epartset,\mu),$ $\nu\left[ W^2 f^2\right] < \infty \}$,
$\consistweight = \{ f \in \mathbb{B}(\epartset), \nu \left[ W^2|f| \right] < \infty \}$, $\gamma(f) =  \nu \left[W^2 f \right]$ for
$f \in \consistweight$, and $\sigma^2(f)= \PVar_\nu \left[ W  f \right]$ for  $f \in \asnormfunc$. The sampling importance resampling
procedure outlined in \cite{rubin:1987} consists in resampling the weighted sample $\{ (\epart{N,i}, W(\epart{N,i}) \}_{1 \leq i \leq M_N}$.
It is recommended that the number $\tilde{M}_N$ of resampled particles be much less than $M_N$ the number of elements in the
importance weighted sample $\{(\epart{N,i},\ewght{N,i} ) \}_{1 \leq i \leq M_N}$ (in \cite[pp. 459]{gelman:rubin:1992} it is suggested  to sample
$\tilde{M}_N = 10$ out of $M_N= 1000$ elements). This suggestion is supported by the results above.
Theorem \ref{thm:LLN-unbiased-sampling} shows that the uniformly weighted sample $\{ (\etpart{N,i},1) \}_{1 \leq i \leq \tilde{M}_N}$
is consistent for $\{\mu,\lone(\epartset,\mu)\}$. Assume that $\beta \geq 1$, Theorem \ref{thm:CLT-independent-sampling} shows that
$\{ (\etpart{N,i},1) \}_{1 \leq i \leq \tilde{M}_N}$ is asymptotically normal for
$(\mu,$ $\tilde{\asnormfunc},$ $\lone(\epartset,\mu)$, $\tilde{\stdnormfunc},$ $\mu,$ $ \sqrt{\tilde{M_N}}\})$
where $\tilde{\asnormfunc} \eqdef \left \{ f \in \lone(\epartset,\mu), \nu( (W + W^2) f^2) < \infty \right\}$ and $\tilde{\stdnormfunc}$ is defined by
$\tilde{\stdnormfunc}^2(f)= \PVar_\mu(f)+  \beta^{-1} \PVar_\nu \left[ W f \right]$.
Suppose now that $\lim_{N \to \infty} M_N/\tilde{M}_N = \infty$ and thus $\beta= \infty$. Then, the limiting variance is the basic Monte Carlo variance
$\PVar_\mu(f)$: the resampled particles can be
thought of as an \iid\ sample from the target distribution $\mu$, providing thus a simple theoretical justification for the SIR procedure.
\end{ex}

The analysis of the deterministic plus residual sampling is more involved.
To carry out the analysis, it is required to specify more precisely the importance weight.
\begin{thm}
\label{thm:resid}
Let $k$ be an integer, $\nu$  be  probability measure on $(\epartset,\epartsigma)$, $\gamma$ be a finite measure on $(\epartset, \epartsigma)$,
$\asnormfunc \subseteq \lone(\epartset,\nu)$, $\consistfunc \subseteq \lone(\epartset,\nu)$ and $\consistweight \subseteq \lone(\Xset,\gamma)$
be proper sets, $\stdnormfunc$ be a non negative function on $\asnormfunc$ and $\funcweight$ be a non negative function on $\epartset$.
Assume that
\begin{enum_i}
\item $\{(\epart{N,i},\funcweight(\epart{N,i}))\}_{1 \leq i \leq M_N}$ is consistent for $(\nu,\consistfunc)$ and
asymptotically normal for $(\nu,$ $\asnormfunc,$ $\consistweight$, $\stdnormfunc,$ $\gamma,$ $\{a_N\})$,
\item $\lim_{N \to \infty} \tilde{M}_N / M_N = \ell$, $\lim_{N \to \infty} a_N= \infty
$ and $\lim_{N \to \infty} a_N^2/\tilde M_N = \beta$, where $\ell \in [0,\infty]$ and $\beta \in [0,\infty]$,
\item \label{item:nechargepas}
$1/\funcweight \in \consistfunc$, and $\nu\left( \ell \nu(1/\funcweight) \funcweight \in \nset \cup \{\infty\} \right)=0$.
\end{enum_i}
Define
\begin{equation}
\label{eq:Wellphi}
W[\ell,\Phi] \eqdef \begin{cases} 1 - \frac{\ell \nu(1/\funcweight) \funcweight}{\lfloor \ell \nu(1/\funcweight) \funcweight \rfloor} & \ell < \infty \\
                                  0 & \ell = \infty \eqsp.
                    \end{cases}
\end{equation}
Then, the following holds true for the uniformly weighted sample $\{(\etpart{N,i},1)\}_{1 \leq i \leq \tilde M_N}$
defined by the deterministic-plus-residual sampling (\eqref{eq:multinomial-sampling-resid-1} and \eqref{eq:multinomial-sampling-resid-2})
\begin{enum_i}
\item If $\beta < 1$, then $\{(\etpart{N,i},1)\}_{1 \leq i \leq \tilde M_N}$ is asymptotically normal for $(\nu,$ $\tilde{\asnormfunc},$ $\consistfunc$, $\tilde{\stdnormfunc},$ $\tilde{\gamma},$ $\{a_N\})$ where
$\tilde{\asnormfunc}$ is given by \eqref{eq:definition-asnormfunc-sampling} and
\begin{align*}
&\tilde{\stdnormfunc}^2(f)\eqdef \beta \nu\left\{W[\ell,\funcweight] \left( f - \frac{\nu\left\{W[\ell,\funcweight] f\right\}}{\nu\{W[\ell,\Phi]\}}\right)^2\right\} +  \stdnormfunc^2\left(f\right) \eqsp, \quad f \in \tilde{\asnormfunc} \eqsp,\\
&\tilde{\gamma}\eqdef \beta \nu \eqsp.
\end{align*}
\item If $\beta \geq 1$, then $\{(\etpart{N,i},1)\}_{1 \leq i \leq \tilde{M}_N}$ is asymptotically normal for  $(\nu,$ $\tilde{\asnormfunc},$ $\consistfunc,$ $\tilde{\stdnormfunc},$ $\tilde{\gamma},$ $\{ \tilde{M}^{1/2}_N \})$
 where $\tilde{\asnormfunc}$ is given by \eqref{eq:definition-asnormfunc-sampling} and
\begin{align*}
&\tilde{\stdnormfunc}^2(f)\eqdef \nu\left\{W[\ell,\funcweight] \left( f - \frac{\nu\left\{W[\ell,\funcweight] f\right\}}{\nu\{W[\ell,\Phi]\}}\right)^2\right\}  +  \beta^{-1} \stdnormfunc^2\left(f\right) \eqsp, \quad f \in \tilde{\asnormfunc} \eqsp,\\
&\tilde{\gamma}\eqdef \nu \eqsp,
\end{align*}
\end{enum_i}
\end{thm}

Because $W[\ell,\Phi] \leq 1$, for any $f \in \tilde{\asnormfunc}$,
\begin{multline*}
\nu\left\{W[\ell,\funcweight] \left( f - \frac{\nu\left\{W[\ell,\funcweight] f\right\}}{\nu\{W[\ell,\Phi]\}}\right)^2\right\} =
\inf_{c \in \rset} \nu \left \{ W[\ell,\funcweight] (f-c)^2 \right\}\\
 \leq \inf_{c \in \rset} \nu \left \{(f-c)^2 \right\} = \PVar_\nu(f) \eqsp,
\end{multline*}
showing that the variance of the residual plus deterministic sampling is always lower than that of the multinomial sampling. These results
extend \cite[Theorem 2]{chopin:2004} that derive an expression of the variance of the residual sampling in a specific case. Note however
the assumption Theorem \ref{thm:resid}-\ref{item:nechargepas} is missing in the statement of \cite[Theorem 2]{chopin:2004}.
This assumption cannot be relaxed, as shown in Section \ref{sec:proof:thm:resid}.

\section{An application to state-space models}
\label{sec:application to state-space models}
In this section, we apply the results developed above to state-space models.
State-space  model  has become a powerful tool for  dynamic systems. In this model,
an underlying state of interest  changes over time and measurements are taken to enable inferences to be made about the state.
The state process  $\{X_k\}_{k \geq 1}$ is a Markov chain on a general state space $(\Xset,\Xsigma)$
with initial distribution $\Xinit$ and  kernel $Q$. The observations $\{ Y_k \}_{k \geq 1}$ are random variables
taking value in a general state space $(\Yset,\Ysigma)$ that are independent conditionally on the state sequence $\{X_k\}_{k \geq 1}$;
in addition, there exists a measure $\lambda$ on $(\Yset,\Ysigma)$, and a transition density function $x \mapsto g(x,y)$, referred to as
the likelihood, such that $\CPP{Y_k \in A}{X_k= x} = \int_A g(x,y) \lambda(dy)$, for all $A \in \Yset$.
The  kernel $Q$ and the likelihood functions $x \mapsto g(x,y)$ are assumed to be known.
These quantities could be time-dependent. Such models have been used in a wide range of applications, including quantitative finance,
engineering and natural sciences and consequently have become of increasing interest to statisticians (see for instance \cite[chapter 3,4]{liu:2001} and \cite{kuensch:2001} for an introduction to that field).

In this paper, we are primarily concerned with the recursive estimation of the (joint) smoothing distribution, \ie\
the conditional distribution of the state sequence $\chunk{X}{1}{k} \eqdef (X_1, \dots, X_k)$ given the $\sigma$-algebra
generated by the observed process from time $1$ to $k$,  \ie\ one is interested in estimating the conditional expectation
\begin{equation}
\label{eq:filteringdistribution}
\post[\Xinit]{k}(\chunk{y}{1}{k},f) \eqdef  \CPE{ f(\chunk{X}{1}{k})}{\chunk{Y}{1}{k}= \chunk{y}{1}{k}} \eqsp, \quad \text{where $f \in \mathbb{B}(\Xset^{k})$} \eqsp.
\end{equation}
We shall consider the case in which the observations have an arbitrary but fixed value $\chunk{y}{1}{k}$, and we drop them from the
notations. We denote $g_k(x)= g(x,y_k)$. A straightforward application of the Bayes formula, shows that, for any $f \in \mathcal{B}^+(\Xset^k)$,
\begin{equation}
\label{eq:joint:smoothing-propagation}
\post[\Xinit]{k}(f) = \frac{\idotsint \post[\Xinit]{k-1}(d \chunk{\tilde{x}}{1}{k-1}) \Q(\tilde{x}_{k-1},d \tilde{x}_k) g_k(\tilde{x}_k) f(\chunk{\tilde{x}}{1}{k})}{\idotsint \post[\Xinit]{k-1}(d \chunk{\tilde{x}}{1}{k-1}) Q(\tilde{x}_{k-1},d \tilde{x}_k) g_k(\tilde{x}_k)} \eqsp.
\end{equation}
In practice, these computations can only be performed in closed-form for linear Gaussian models and for
finite state-space models.
Many approximations schemes have been proposed in the literature to tackle this problem, but most of the known solutions
either suffer from poor performance and / or instability (Extended Kalman Filter, Gaussian sum filter, etc.),
or are prohibitively expensive to implement (grid-based solutions).  In the Monte-Carlo framework, we approximate the posterior distribution $\post[\Xinit]{k}$
at each iteration $k$ by means of a weighted sample $\{ (\epart[k]{N,i}, \ewght[k]{N,i}) \}_{1 \leq k \leq M_N}$, where
the superscript $k$ indicates the iteration index.  Note that $\post[\Xinit]{k}$ is defined on $(\Xset^{k+1}, \Xsigma[k+1])$ and thus that the
points $\epart[k]{N,i}$  belong to $\Xset^{k+1}$ (these are often referred to as \emph{path particle} in the literature; see \cite{delmoral:2004}).

To apply the results presented in section \ref{sec:mainresults}, it is first required to define a transition kernel $\LTAR_{k-1}$ satisfying
\eqref{eq:definition-mu} with $\nu= \post[\Xinit]{k-1}$, $(\epartset,\epartsigma)= (\Xset^k, \Xsigma[k])$,
$\mu= \post[\Xinit]{k}$ and $(\etpartset,\etpartsigma) = (\Xset^{k+1},\Xsigma[k+1])$, \ie\ for any $f \in \mathbb{B}^+(\Xset^{k+1})$,
\begin{equation}
\label{eq:definition-flow}
\post[\Xinit]{k}(f) = \frac{\idotsint \post[\Xinit]{k-1}(d \chunk{x}{1}{k-1}) \LTAR_{k-1}(\chunk{x}{1}{k-1},d \chunk{\tilde{x}}{1}{k})
f(\chunk{\tilde{x}}{1}{k})}{\idotsint \post[\Xinit]{k-1}(d \chunk{x}{1}{k-1}) \LTAR_{k-1}(d \chunk{x}{1}{k-1},\Xset^{k+1})}
\end{equation}
In the second step, we must choose a proposal kernel $\KISS_{k-1}$ satisfying
%The instrumental distribution $\XinitIS$ for the initial state dominates the filtering distribution $\filt[\Xinit]{1}$, $\filt[\Xinit]{1} \ll \XinitIS$.
\begin{equation}
\label{eq:dominationkernel}
 \LTAR_{k-1}(\chunk{x}{1}{k-1},\cdot) \ll \KISS_{k-1}(\chunk{x}{1}{k-1},\cdot) \eqsp, \quad \text{for any  $\chunk{x}{1}{k-1} \in \Xset^{k}$.}
\end{equation}
There are many possible choices, which will be associated with different algorithms proposed in the literature. The first obvious
choice consists in setting, for any $f \in \mathcal{B}(\Xset^k)$,
\begin{equation}
\label{eq:definitionKUN}
\LTAR_{k-1}(\chunk{x}{1}{k-1},f)= \KUN_{k-1}(\chunk{x}{1}{k-1},f)=  \int \Q(x_{k-1},d \tilde{x}_k) g_k(\tilde{x}_k) f(\chunk{x}{1}{k-1},\tilde{x}_k)  \eqsp.
\end{equation}
Note that, by construction, the kernel $\KUN_{k-1}(\chunk{x}{1}{k-1},\cdot)$ leaves the coordinates $\chunk{x}{1}{k-1}$.
 This corresponds to the so-called sequential importance sampling algorithm. The first obvious choice is that of setting $\KISS_{k-1} = \KPRI_{k-1}$, where the kernel $\KPRI_{k-1}$ is defined, for any $f \in \mathbb{B}^+(\Xset^{k})$,
\begin{equation}
\label{eq:definition-KPRI}
\KPRI_{k-1}(\chunk{x}{1}{k-1},f) = \int \Q(x_{k-1},d \tilde{x}_k)f(\chunk{x}{1}{k-1},\tilde{x}_k) \eqsp.
\end{equation}
With this particular choice, for any $\chunk{x}{1}{k-1} \in \Xset^k$ and $\chunk{\tilde{x}}{1}{k} \in \Xset^k$,
\begin{equation}
\label{eq:radon-nykodim-transition-a-priori}
\frac{d \KUN_{k-1}(\chunk{x}{1}{k-1},\cdot)}{d \KPRI_{k-1}(\chunk{x}{1}{k-1},\cdot)}(\chunk{\tilde{x}}{1}{k})   \propto g_k(\tilde{x}_k) \eqsp.
\end{equation}
Note that the incremental weight $g_k(\tilde{x}_k)$ \emph{does not depend on} $\chunk{x}{1}{k-1}$, that is,
on the past path particle. The use of the prior kernel $\KISS_{k-1}= \KPRI_{k-1}$
is popular because sampling from the prior kernel $\KPRI_{k-1}$ is often straightforward,
and computing the incremental weight simply amounts to evaluating the
conditional likelihood of the new observation given the current particle
position. Often, significant gain can be expected by taking into account the new observation in the
mutation kernel $\KISS_{k-1}$  (see for instance \cite{liu:chen:1998}, \cite{doucet:godsill:andrieu:2000} and
\cite{doucet:defreitas:gordon:2001}). A possible choice is to set the instrumental kernel $\KISS_{k-1}(\chunk{x}{1}{k-1},\cdot)$ as the conditional
distribution of the state $X_{k}$ given the previous state and the current observation, \ie\ to set $\KISS_{k-1}= \KOPT_{k-1}$, where
for any $f \in \mathbb{B}^+(\Xset^{k})$ and  $\chunk{x}{1}{k-1} \in \Xset^{k-1}$,
$$
\KOPT_{k-1}(\chunk{x}{1}{k-1},f) = \frac{\int \Q(x_{k-1},d \tilde{x}_k) g_k(\tilde{x}_k) f(\chunk{x}{1}{k-1},\tilde{x}_k)}{\int \Q(x_{k-1},dx_k) g_k(x_k)}   \eqsp,
$$
which is often termed the \emph{optimal kernel}. The incremental importance weight is thus given by
\begin{equation}
\label{eq:radon-nykodim-transition-optimal}
\frac{d \KUN_{k-1}(\chunk{x}{1}{k-1},\cdot)}{d \KOPT_{k-1}(\chunk{x}{1}{k-1},\cdot)}(\chunk{\tilde{x}}{1}{k}) = \int \Q(x_{k-1},dx_k) g_k(x_k) \eqsp ,
\end{equation}
which depends on the last coordinate of the past value of the path particle $\chunk{x}{1}{k-1}$ but not on the current value of the particle
offspring $\chunk{\tilde{x}}{1}{k}$.
Most often, it is not easy to sample from $\KOPT_{k-1}$ or to compute the importance weights $\KUN_{k-1}(\epart[k-1]{N,j},\Xset)$.
A possible solution consists in sampling from $\KOPT_{k-1}$ by accept-reject  (see \cite{tanizaki:1999} and \cite{kuensch:2003});
another solution is to  approximate  the optimal kernel  by using more or less sophisticated tricks
(see for instance \cite{shephard:pitt:1997}, \cite{doucet:godsill:andrieu:2000}, \cite{tanizaki:2001}).

Other choices for $\LTAR_{k-1}$ and $\KISS_{k-1}$ are possible. For example, \cite{gilks:berzuini:2001b} have introduced a variant
of this algorithm, the \emph{resample-move} algorithm, in which the whole path of the particles are mutated. This technique allows to
combat the progressive impoverishment of the system of particles as the dynamic process evolves.
The construction goes as follows. Let $P_{k-1}$ be a kernel on $(\Xset^{k}, \Xsigma[k])$ such that $\post[\Xinit]{k-1}$ is an invariant distribution for $P_{k-1}$, \ie\
$\post[\Xinit]{k-1} P_{k-1}= \post[\Xinit]{k-1}$. Such kernel can be constructed easily by using for example the Metropolis-Hastings
construction (which automatically guarantees the detailed balance condition). Then, define the kernel $\LTAR_{k-1}$, for any $f \in \mathcal{B}(\Xset^{k})$ and
$\chunk{x}{1}{k-1} \in \Xset^{k-1}$, by
\begin{equation}
\label{eq:newdefinitionkernel}
\LTAR_{k-1}(\chunk{x}{1}{k-1},f)= \idotsint P_{k-1}(\chunk{x}{1}{k-1},d \chunk{\tilde{x}}{1}{k-1}) Q(\tilde{x}_{k-1},d \tilde{x}_k) g_k(\tilde{x}_k) f(\chunk{\tilde{x}}{1}{k}) \eqsp.
\end{equation}
Because $P_{k-1}$ is invariant for $\post[\Xinit]{k-1}$, for any $f \in \mathcal{B}(\Xset^{k})$,
\begin{multline*}
\idotsint \post[\Xinit]{k-1}( d \chunk{x}{1}{k-1}) P_{k-1}( \chunk{x}{1}{k-1}, d \chunk{\tilde{x}}{1}{k-1}) Q(\tilde{x}_{k-1},d \tilde{x}_k) g_k(\tilde{x}_k) f(\chunk{\tilde{x}}{1}{k}) = \\
\idotsint \post[\Xinit]{k-1}( d \chunk{\tilde{x}}{1}{k-1})  Q(\tilde{x}_{k-1},d \tilde{x}_k) g_k(\tilde{x}_k) f(\chunk{\tilde{x}}{1}{k}) \eqsp,
\end{multline*}
showing that \eqref{eq:definition-flow}. With this definition of $\LTAR_{k-1}$, the condition \eqref{eq:dominationkernel} is satisfied for
example by the kernel $\KISS_{k-1}$ given, for any $f \in \mathcal{B}(\Xset^{k})$,
\[
\KISS_{k-1}(\chunk{x}{1}{k-1},f)= \idotsint P_{k-1}(\chunk{x}{1}{k-1},d \chunk{\tilde{x}}{1}{k-1}) Q(\tilde{x}_{k-1},d \tilde{x}_k) f(\chunk{\tilde{x}}{1}{k}) \eqsp.
\]
Other possible choices and implementations issues are given in \cite{gilks:berzuini:2001b}.

We proceed from the weighted sample $\{ (\epart[k-1]{N,i}, \ewght[k-1]{N,i}) \}_{1 \leq i \leq M_N}$ targeting $\post[\Xinit]{k-1}$ to $\{(\epart[k]{N,i}, \ewght[k]{N,i}) \}_{1 \leq i \leq M_N}$ targeting $\post[\Xinit]{k}$ as follows.
To keep the discussion simple, it is  assumed that each particle gives birth to a single offspring.
In the mutation step, we draw $\{ \etpart[k]{N,i} \}_{1 \leq i \leq M_N}$ conditionally independently given $\mcf[k-1]{N}$ with distribution
given, for any $f \in \mathbb{B}^+(\Xset^{k})$ by
\begin{equation}
\label{eq:mutation-step:mutation/selection}
\CPE{f(\etpart[k]{N,i})}{\mcf[k-1]{N}} = \KISS_{k-1}(\epart[k-1]{N,i},f) = \idotsint \KISS_{k-1}(\epart[k-1]{N,i},d \chunk{x}{1}{k}) f(\chunk{x}{1}{k}) \eqsp,
\end{equation}
where  $i=1, \dots, M_N$. Next we assign to the particle $\etpart[k]{N,i}$, $i=1, \dots, M_N$, the importance weight
\begin{equation}
\label{eq:definitionWk}
\etwght[k]{N,i}=  \ewght[k-1]{N,i} W_{k-1}(\epart[k-1]{N,i},\etpart[k]{N,i}) \ \  \text{with} \ \
W_{k-1}(\chunk{x}{1}{k-1},\chunk{\tilde{x}}{1}{k}) = \frac{d \LTAR_{k-1}(\chunk{x}{1}{k-1},\cdot)}{d \KISS_{k-1}(\chunk{x}{1}{k-1},\cdot)}(\chunk{\tilde{x}}{1}{k})\eqsp.
\end{equation}
Instead of resampling at each iteration (which is the assumption upon which most of the asymptotic analysis have been carried out so far),
we rejuvenate the particle system only when the importance weights are too skewed. As discussed in \cite[section 4]{kong:liu:wong:1994}, a sensible approach
is to try to monitor the coefficient of variations of weight, defined by
\[
\left[ \CV{N}{k} \right]^2 \eqdef \frac{1}{M_N} \sum_{i=1}^{M_N} \left( \frac{M_N \etwght[k]{N,i}}{\tsumweight[k]{N}} -  1\right)^2 \eqsp.
\]
The coefficient of variation is minimal when the normalized importance weights $\etwght[k]{N,i}/ \tsumweight[k]{N}$, $i=1, \dots, M_N$,
are all equal to $1/M_N$, in which case $\CV{N}{k} = 0$. The maximal value of $\CV{N}{k}$ is
$\sqrt{M_N-1}$, which corresponds to one of the normalized weights being one and
all others being null. Therefore, the coefficient of variation is often
interpreted as a measure of the number of ineffective particles (those that do
not significantly contribute to the estimate). A related criterion with a
simpler interpretation is the so-called \emph{effective sample  size}\index{Effective sample size} $\mathrm{ESS}_{N}^{(k)}$ \citep{liu:1996},
defined as
\begin{equation}\label{eq:EffectiveSampleSize}
\mathrm{ESS}_{N}^{(k)} = \left[ \sum_{i=1}^{M_N} \left(\frac{\etwght[k]{N,i}}{\tsumweight[k]{N}} \right)^2 \right]^{-1} \eqsp ,
\end{equation}
which varies between 1 (all weights null but one) and $M_N$ (equal weights).
It is straightforward to check the relation
$$
\mathrm{ESS}_{N}^{(k)} = \frac{M_N}{1+\left[ \CV{N}{k} \right]^2} \eqsp .
$$
The effective sample size may be understood as a proxy for the equivalent number of \iid\ samples at time $k$.
Some additional insights and heuristics about the coefficient of variation are given by \cite{liu:chen:1995}.
If the coefficient of variation of the importance weights (or equivalently, if the ratio of the effective sample size to the total
number of particles, $\mathrm{ESS}_{N}^{(k)}/M_N$) crosses a threshold we rejuvenate the particle system.  More precisely if at time $k$,
$\CV{N}{k} \geq \kappa$ (such time index are called \emph{dynamic check points} in \cite{liu:chen:logvinenko:2001}),
we draw $I_{k}^{N,1}, \dots, I_{k}^{N,M_N}$ conditionally independently given $\mctf[k]{N}= $ $\mcf[k-1]{N} \vee$
$\sigma \bigl( \{ \etpart[k]{N,i},\etwght[k]{N,i} \}_{1 \leq i \leq M_N} \bigr)$,
with distribution
\begin{equation}
\label{eq:selection-step:mutation/selection}
\CPP{I_{k}^{N,i}=j}{\mctf[k]{N}} =  \etwght[k]{N,j} / \tsumweight[k]{N} \eqsp, \quad i= 1, \dots, M_N \eqsp, j=1, \dots, M_N
\end{equation}
and we set
\begin{equation}
\label{eq:afterresampling}
\epart[k]{N,i} = \etpart[k]{N,I_{k}^{N,i}} \quad \text{and} \quad \ewght[k]{N,i}= 1 \eqsp, i=1, \dots, M_N \eqsp.
\end{equation}
If $\CV{N}{k} < \kappa$, we simply copy the mutated path particles: $(\epart[k]{N,i}, \ewght[k]{N,i})= (\etpart[k]{N,i},\etwght[k]{N,i})$, $i=1, \dots, M_N$.
In both cases, we set $\mcf[k]{N} = \mctf[k]{N} \vee $ $\sigma( \{ (\epart[k]{N,i}, \ewght[k]{N,i}) \}_{1 \leq i \leq M_N}$.
We consider here only multinomial resampling, but the deterministic plus residual sampling
can be applied as well.

%\begin{assumSMC}
%\label{assumSMC:bootstrapI}
%$\Xinit(g_0)> 0$, $\sup_{\x\in\Xset} \g_k(x) < \infty$ for all $k\geq 0$ and
%\begin{equation}
%\label{item:bootstrapI:assum-2}
%\int_\Xset \Q(x,dx') g_k(x') > 0 \eqsp,  \text{for all $\x \in \Xset$ and  $k\geq 0$} \eqsp.
%\end{equation}
%\end{assumSMC}

\begin{thm}
\label{thm:bootstrapfilterI}

For any $k > 0$, let $\LTAR_{k}$ and $\KISS_k$ be  transition kernels from $(\Xset^k,\Xsigma[k])$ to $(\Xset^{k+1},\Xsigma[k+1])$ satisfying \eqref{eq:definition-flow} and \eqref{eq:dominationkernel}, respectively. Assume that the equally weighted sample $\{( \epart[1]{N,i},1) \}_{1 \leq i \leq M_N}$ is consistent for $\{\filt[\Xinit]{1},\lone \left(\Xset,\filt[\Xinit]{1}\right)\}$ and asymptotically normal for
$(\filt[\Xinit]{1},$ $\asnormfunc_1,$ $\consistweight_1,$ $\stdnormfunc_1,$ $\filt[\Xinit]{1},$ $\{M_N^{1/2}\})$ where $\asnormfunc_1$ and $\consistweight_1$ are proper sets and define recursively $(\asnormfunc_k)$ and $(\consistweight_k)$ by
\begin{align*}
&
\asnormfunc_k \eqdef \left\{ f \in \ltwo \left( \Xset^{k}, \post[\Xinit]{k} \right), \quad  \LTAR_{k-1}(\cdot,f) \in \asnormfunc_{k-1} \eqsp,
\KISS_{k-1}(\cdot, W_{k-1}^2 f^2) \in \consistweight_{k-1} \right\} \eqsp, \\
&
\consistweight_k \eqdef \left\{ f \in \lone \left(\Xset^{k},\post[\Xinit]{k} \right), \quad \KISS_{k-1}( \cdot, W_{k-1}^2|f|) \in \consistweight_{k-1} \right\} \eqsp.
\end{align*}
Assume in addition that for any $k\geq 1$,  $\KISS_{k}( \cdot, W_{k}^2) \in \consistweight_{k}$. Then
for any $k\geq 1$, $(\asnormfunc_k)$ and $(\consistweight_k)$ are proper sets and $\{ (\epart[k]{N,i},\ewght[k]{N,i}) \}_{1 \leq i \leq M_N}$  is consistent for $\{\filt[\Xinit]{k},\lone \left(\Xset,\filt[\Xinit]{k}\right)\}$ and asymptotically normal for $(\filt[\Xinit]{k},$
$\asnormfunc_k,$
$\consistweight_k,$ $\stdnormfunc_k,$ $\gamma_k,$ $\{M_N^{1/2}\})$,
where the functions  $\stdnormfunc_k$ and the measure $\gamma_k$ are given by
\begin{align*}
& \stdnormfunc^2_{k}(f) = \varepsilon_k \PVar_{\post[\Xinit]{k}}(f)  \\
& \quad \quad +\frac{\stdnormfunc^2_{k-1}(\LTAR_{k-1} \{ f - \post[\Xinit]{k}(f) \} ) +  \gamma_{k-1} \KISS_{k-1} \left[ \left\{ W_{k-1}  f - \KISS_{k-1}(\cdot, W_{k-1} f) \right\}^2 \right] }
{\left\{  \post[\Xinit]{k-1}\LTAR_{k-1}(\Xset^{k})\right\}^2} \\
&\gamma_k(f)= \varepsilon_k \post[\Xinit]{k} + (1- \varepsilon_k) \frac{\gamma_{k-1} \KISS_{k-1}(W_{k-1}^2 f)}{[\post[\Xinit]{k-1} \LTAR_{k-1}(\Xset^k) ]^2}
\end{align*}
where $W_k$ is defined in \eqref{eq:definitionWk} and
\[
\varepsilon_k \eqdef \1 \left\{ [\post[\Xinit]{k-1} \LTAR_{k-1}(\Xset^k)]^{-2} \gamma_{k-1} \KISS_{k-1}(W_{k-1}^2) \geq 1 + \kappa^2 \right\}
\]
\end{thm}

\begin{proof}
Recall that the algorithm proceeds as follows. The weighted sample $\{ (\epart[k-1]{N,i}, \ewght[k-1]{N,i}) \}_{1 \leq i \leq M_N}$ are mutated into $\{ (\etpart[k]{N,i}, \etwght[k]{N,i}) \}_{1 \leq i \leq M_N}$  as described by (\ref{eq:mutation-step:mutation/selection}) and (\ref{eq:definitionWk}) (see also section \ref{subsec:mutation}). The resulting particles are then resampled according to the multinomial resampling algorithm (as described in Section \ref{subsec:resampling}) so as to obtain a family of equally weighted particles that we denote by $\{ ( \ebpart[k]{N,i},1) \}_{1 \leq i \leq M_N}$.
$$
 (\epart[k-1]{N,i}, \ewght[k-1]{N,i}) \overset{\mbox{\tiny mutation}}{\longrightarrow}(\etpart[k]{N,i}, \etwght[k]{N,i}) \overset{\mbox{\tiny resampling}}{\longrightarrow}( \ebpart[k]{N,i},1)
$$
We then assign
\begin{equation*}
  (\epart[k]{N,i}, \ewght[k]{N,i})=
\begin{cases}
 (\etpart[k]{N,i}, \etwght[k]{N,i})  & \mbox{if}\ \CV{N}{k}< \kappa\\
 ( \ebpart[k]{N,i},1)   & \mbox{otherwise}
\end{cases}
\end{equation*}

The proof now follows by induction. Assume that for some $k>1$, the weighted sample $\{ (\epart[k-1]{N,i}, \ewght[k-1]{N,i}) \}_{1 \leq i \leq M_N}$ is consistent for  $\{\filt[\Xinit]{k-1},\lone \left(\Xset,\filt[\Xinit]{k-1}\right)\}$ and asymptotically normal for  $(\filt[\Xinit]{k-1},$
$\asnormfunc_{k-1},$
$\consistweight_{k-1},$ $\stdnormfunc_{k-1},$ $\gamma_{k-1},$ $\{M_N^{1/2}\})$. By Theorem \ref{thm:LLN-mutation} and \ref{thm:CLT-mutation},  $(\etpart[k]{N,i}, \etwght[k]{N,i})_{1 \leq i \leq M_N}$ is consistent for  $\{\filt[\Xinit]{k},\lone \left(\Xset,\filt[\Xinit]{k}\right)\}$ and asymptotically normal for $(\filt[\Xinit]{k},$ $\tilde \asnormfunc_k,$ $\tilde \consistweight_k,$ $\tilde \stdnormfunc_k,$ $\tilde \gamma_k,$ $\{M_N^{1/2}\})$ where  $\tilde \asnormfunc_k,$ $\tilde \consistweight_k,$ $\tilde \stdnormfunc_k,$ $\tilde \gamma_k,$ are defined from $\asnormfunc_k,$ $\consistweight_k,$ $\stdnormfunc_k,$ $\gamma_k,$ using Theorem \ref{thm:CLT-mutation}. And by Theorem \ref{thm:LLN-unbiased-sampling} and \ref{thm:CLT-independent-sampling},  $( \ebpart[k]{N,i},1)_{1 \leq i \leq M_N}$ is  consistent for $\{\filt[\Xinit]{k},\lone \left(\Xset,\filt[\Xinit]{k}\right)\}$ and  asymptotically normal for $(\filt[\Xinit]{k},$ $\bar \asnormfunc_k,$ $\bar \consistweight_k,$ $\bar \stdnormfunc_k,$ $\bar \gamma_k,$ $\{M_N^{1/2}\})$ where  $\bar \asnormfunc_k,$ $\bar \consistweight_k,$ $\bar \stdnormfunc_k,$ $\bar \gamma_k,$ are defined from $\tilde \asnormfunc_k,$ $\tilde \consistweight_k,$ $\tilde \stdnormfunc_k,$ $\tilde \gamma_k,$ using Theorem \ref{thm:CLT-independent-sampling}. The asymptotic normality of  $(\etpart[k]{N,i}, \etwght[k]{N,i})_{1 \leq i \leq M_N}$ and   $( \ebpart[k]{N,i},1)_{1 \leq i \leq M_N}$, combined with:
$$
\left[\CV{N}{k}\right]^2=M_N  \sum_{i=1}^{M_N} \left(\frac{\etwght[k]{N,i}}{\tsumweight[k]{N}} \right)^2 -1 \plim \tilde \gamma_k(1)-1\quad \mbox{and} \quad  \epsilon_k=\1\left\{\tilde \gamma_k(1)-1>\kappa^2 \right\} \eqsp,
$$
complete the proof.

\end{proof}

\appendix
\section{Weak Limits Theorems for Triangular Array}
\label{sec:WeakLimitsTriangularArray}

This section summarizes various  limit theorems for triangular arrays of dependent random variables
To keep the technical assumptions in our main theorems minimal, we derive these limit theorems
under assumptions that are weaker than the one typically used in the literature (see \cite{hall:heyde:1980}).
Through our general results can be obtained by weakening the assumptions and adapting the proofs for
triangular array of dependent random variables given in \cite{dvoretzky:1972} and \cite{mcleish:1974}, we prefer to develop them here independently.
We hope that the greater accessibility of the proofs will compensate for this sacrifice of brevity.
Let $(\Omega,\mcf{},\PP)$ be a probability space, let $X$ be a random variable and let $\mcg{}$ be a a sub-$\sigma$ field of $\mathcal{F}$.
Define $X^+ \eqdef \max(X,0)$ and $X^- \eqdef -\min(X,0)$. Following \cite[Section II.7]{shiryaev:1996}, if
\[
\min( \CPE{X^+}{\mcg{}}, \CPE{X^-}{\mcg{}} ) < \infty \eqsp \text{, \PP-\as}
\]
(a version of) the conditional expectation of $X$ given $\mcg{}$
is defined by
$$
\CPE{X}{\mcg{}}= \CPE{X^+}{\mcg{}} - \CPE{X^-}{\mcg{}} \eqsp,
$$
where, on the $\PP$-null-set of sample points for which $\CPE{X^+}{\mcg{}}= \CPE{X^-}{\mcg{}}= \infty$, the difference $\CPE{X^+}{\mcg{}}
- \CPE{X^-}{\mcg{}}$ is given an arbitrary value, for
instance, zero. In particular, if $\CPE{|X|}{\mcg{}} < \infty$ $\PP$-\as\ then $\CPE{X^+}{\mcg{}} < \infty$ and
$\CPE{X^-}{\mcg{}} < \infty$ $\PP$-\as\ and we may always define the conditional expectation in this context.

Let $\{M_N \}_{N \geq 0}$ be a sequence of positive integers
satisfying $\lim_{N \to \infty} M_N = \infty$. Without loss of generality, we will assume that $\{M_N \}_{N \geq 0}$ is non decreasing. Let $\{U_{N,i} \}_{1
\leq i \leq M_N}$  be a triangular array of random variables on
$(\Omega, \mcf{}, \PP)$. Let $\{ \mcf{N,i} \}_{0 \leq i \leq M_N}$
be a triangular array of of sub-sigma-fields of $\mcf{}$ such that
for each $N$ and each $i=1, \dots, M_N$, $U_{N,i}$ is
$\mcf{N,i}$-measurable and $\mcf{N,i-1} \subseteq \mcf{N,i}$.

\begin{prop}
\label{prop:LLN-NZM} Assume that $\CPE{|U_{N,j}|}{\mcf{N,j-1}} < \infty$ $\PP$-\as\ for any $N$ and any $j=1, \dots, M_N$, and
\begin{align}
\label{eq:LLN-NZM:tightness}
&\sup_{N} \PP \left( \sum_{j=1}^{M_N} \CPE{|U_{N,j}|}{\mcf{N,j-1}} \geq \lambda \right) \to 0 && \text{ as $\lambda \to \infty$}\\
\label{eq:LLN-NZM:negligeability}
&\sum_{j=1}^{M_N} \CPE{|U_{N,j}| \1 {\{ |U_{N,j}| \geq \epsilon \}}}{\mcf{N,j-1}} \plim 0 &&
\text{for any positive $\epsilon$} \eqsp.
\end{align}
Then,
$
\max_{1 \leq i \leq M_N}\left|\sum_{j=1}^{i} U_{N,j} -  \sum_{j=1}^{i} \CPE{U_{N,j}}{\mcf{N,j-1}}\right|  \plim 0 \eqsp.
$
\end{prop}

\begin{proof}
Assume first that for each $N$ and each $i=1, \dots, M_N$, $U_{N,i} \geq 0$, $\PP$-\as. By \cite[Lemma 3.5]{dvoretzky:1972}, we have that for any constants $\epsilon$ and $\eta > 0$,
\[
\PP \left[ \max_{1 \leq i \leq M_N} U_{N,i} \geq \epsilon \right] \leq \eta + \PP \left[ \sum_{i=1}^{M_N} \CPP{U_{N,i} \geq \epsilon}{\mcf{N,i-1}} \geq \eta \right]\eqsp.
\]
From the conditional version of the Chebyshev identity,
\begin{equation} \label{eq:controleDuMax}
\PP \left[ \max_{1 \leq i \leq M_N} U_{N,i} \geq \epsilon \right] \leq \eta + \PP \left[ \sum_{i=1}^{M_N} \CPE{U_{N,i} \1 \{ U_{N,i} \geq \epsilon \}}{\mcf{N,i-1}} \geq \eta \epsilon \right] \eqsp.
\end{equation}
Let $\epsilon$ and $\lambda >0$ and define $\bar{U}_{N,i}$ by
\[
\bar{U}_{N,i} \eqdef U_{N,i} \1 \{ U_{N,i} < \epsilon \} \, \1 \left \{ \sum_{j=1}^{i} \CPE{U_{N,j}}{\mcf{N,j-1}} < \lambda \right \} \eqsp.
\]
For any $\delta > 0$,
\begin{align*}
&\PP \left( \max_{1 \leq i \leq M_N} \left| \sum_{j=1}^i U_{N,j} - \sum_{j=1}^i \CPE{U_{N,j}}{\mcf{N,j-1}} \right| \geq 2 \delta \right)\\
& \qquad \leq \PP \left( \max_{1 \leq i \leq M_N} \left| \sum_{j=1}^i \bar{U}_{N,j} - \sum_{j=1}^i \CPE{\bar{U}_{N,j}}{\mcf{N,j-1}} \right| \geq  \delta \right) \\
& \qquad + \PP \left ( \max_{1 \leq i \leq M_N} \left| \sum_{j=1}^{i}  U_{N,j} - \bar{U}_{N,j} - \sum_{j=1}^{i}  \CPE{U_{N,j} - \bar{U}_{N,j}}{\mcf{N,j-1}} \right| \geq \delta\right) \eqsp.
\end{align*}
The second term in the \rhs\ is bounded by
\begin{multline*}
\PP \left( \max_{1 \leq i \leq M_N} U_{N,i} \geq \epsilon \right) + \PP \left( \sum_{j=1}^{M_N} \CPE{U_{N,j}}{\mcf{N,j-1}} \geq \lambda \right) + \\
\PP \left( \sum_{j=1}^{M_N} \CPE{U_{N,j} \1 \{ U_{N,j} \geq \epsilon \}}{\mcf{N,j-1}} \geq  \delta \right)
\end{multline*}
Eqs. \eqref{eq:LLN-NZM:negligeability} and \eqref{eq:controleDuMax} imply that the first and last terms in the last expression
converge to zero for any $\epsilon>0$ and \eqref{eq:LLN-NZM:tightness} implies that the second term may be arbitrarily small by choosing for $\lambda$ sufficiently large. Now, by the Doob maximal inequality,
\begin{multline*}
\PP \left( \max_{1 \leq i \leq M_N} \left| \sum_{j=1}^{i} \bar{U}_{N,j} - \CPE{\bar{U}_{N,j}}{\mcf{N,j-1}} \right| \geq \delta \right)\\
 \leq
\delta^{-2} \sum_{j=1}^{M_N} \PE \left( \bar{U}_{N,j} - \CPE{\bar{U}_{N,j}}{\mcf{N,j-1}} \right)^2
\end{multline*}
This last term does not exceed
%\begin{multline*}
$$
\delta^{-2} \sum_{i=1}^{M_N} \PE ( \bar{U}_{N,j}^2 ) \leq \delta^{-2} \epsilon \sum_{j=1}^{M_N} \PE[ \bar{U}_{N,j}] \leq
\delta^{-2} \epsilon \PE \left[ \sum_{j=1}^{M_N} \CPE{\bar{U}_{N,j}}{\mcf{N,j-1}} \right] \leq \delta^{-2} \epsilon \lambda \eqsp.
$$
%\end{multline*}
Since $\epsilon$ is arbitrary, the proof follows for $U_{N,j} \geq 0$, $\PP$-\as, for each $N$ and $j=1, \dots, M_N$. The proof extends to
an arbitrary triangular array $\{ U_{N,j} \}_{1 \leq i \leq M_N}$ by applying the preceding result to $\{ U^+_{N,j} \}_{1 \leq i \leq M_N}$ and
$\{U^-_{N,j} \}_{1 \leq j \leq M_N}$.
\end{proof}

\begin{lem}
\label{lem:CLT:TriangularArrayLindeberg}
Assume that for all $N$ and   $i=1,\dots,M_N$, $\sigma_{N,i}^2 \eqdef \CPE{U_{N,i}^2}{\mcf{N,i-1}} < \infty$ $\PP$-\as\ and
$\CPE{U_{N,i}}{\mcf{N,i-1}}=0$, and for all $\epsilon>0$,
\begin{equation}
\label{eq:lem:lindeberg-condition}
\sum_{i=1}^{M_N} \CPE{U_{N,i}^2 \1 { \{ |U_{N,i}| \geq \epsilon  \}}}{\mcf{N,0}} \plim 0 \eqsp.
\end{equation}
Then, if any of the two following conditions holds
\begin{enumerate}
\item [(i)]  $\sum_{i=1}^{M_N} \sigma_{N,i}^2 =1$,
\item [(ii)] $\sigma_{N,i}$ is $\mcf{N,0}$-measurable and $\{\sum_{j=1}^{M_N} \sigma_{N,j}^2\}_{N \geq 0}$ is tight,
\end{enumerate}
 then for any real $u$,
\begin{equation}
\label{eq:limittcl}
\CPE{\exp \left(  \rmi u \sum_{j=1}^{M_N} U_{N,j} \right)}{\mcf{N,0}} - \exp \left( - (u^2/2) \sum_{j=1}^{M_N} \sigma_{N,j}^2 \right) \plim 0 \eqsp.
\end{equation}
\end{lem}
\begin{proof}
Write the following decomposition (with the convention $\sum_{j=a}^b=0$ if $a>b$):
\begin{equation*}
  \mathrm{e}^{ \rmi u \sum_{j=1}^{M_N} U_{N,j}} - \mathrm{e}^{-\frac{u^2}{2} \sum_{j=1}^{M_N} \sigma_{N,j}^2}= \sum_{l =1}^{M_N} \mathrm{e}^{ \rmi u \sum_{j=1}^{l-1} U_{N,j}}\left(\mathrm{e}^{\rmi u U_{N,l}}- \mathrm{e}^{-\frac{u^2}{2} \sigma_{N,l}^2} \right)\mathrm{e}^{-\frac{u^2}{2}\sum_{j=l+1}^{M_N} \sigma_{N,j}^2}\eqsp.
\end{equation*}
Moreover, if (i) or (ii) holds,  $\sum_{j=1}^{l-1} U_{N,j}$ and  $\sum_{j=l+1}^{M_N} \sigma_{N,j}^2$ are $\mcf{N,l-1}$-measurable random variables. To see this, write if (i) holds,  $\sum_{j=l+1}^{M_N} \sigma_{N,j}^2=1-\sum_{j=1}^{l} \sigma_{N,j}^2$ and if (ii) holds, note that $\sum_{j=l+1}^{M_N} \sigma_{N,j}^2$ is a $\mcf{N,0}$-measurable random variable and $\mcf{N,0}\subset \mcf{N,l-1}$ .  This implies that
\begin{multline}
\label{eq:lem:CLT:rel1}
\left| \CPE{\exp \left( \rmi u \sum_{j=1}^{M_N} U_{N,j} \right) - \exp\left( - (u^2/2)\sum_{j=1}^{M_N} \sigma_{N,j}^2  \right)}{\mcf{N,0}}\right| \\
\leq \sum_{l=1}^{M_N} \CPE{\left| \CPE{\exp(\rmi u U_{N,l})}{\mcf{N,l-1}} - \exp(-u^2 \sigma_{N,l}^2/2 )\right|}{\mcf{N,0}} \eqsp.
\end{multline}
For any $\epsilon > 0$,
\begin{align*}
&\left| \CPE{\exp \left( \rmi u  U_{N,l}\right) }{\mcf{N,l-1}} - 1 - \frac{1}{2} u^2 \sigma_{N,l}^2  \right|\\ \nonumber
& \quad \leq \frac{1}{6} |u|^3 \CPE{|U_{N,l}|^3 \1 \{ |U_{N,l}| \leq \epsilon \}}{\mcf{N,l-1}} + u^2 \CPE{U_{N,l}^2 \1 \{ |U_{N,l}| > \epsilon \}}{\mcf{N,l-1}} \\
& \quad \leq \frac{1}{6} \epsilon |u|^3 \sigma_{N,l}^2 + u^2 \CPE{U_{N,l}^2 \1 \{ |U_{N,l}| > \epsilon \}}{\mcf{N,l-1}} \eqsp, \nonumber
\end{align*}
and thus,
\begin{multline} \label{eq:lem:CLT:rel2Bis}
\CPE{ \sum_{l=1}^{M_N} \left| \CPE{\exp \left( \rmi u  U_{N,l} \right) - 1 - \frac{1}{2} u^2 \sigma_{N,l}^2 }{\mcf{N,l-1}}\right|}{\mcf{N,0}} \\
\leq \frac{1}{6} \epsilon |u|^3\sum_{l=1}^{M_N} \sigma_{N,l}^2 + u^2 \sum_{l=1}^{M_N} \CPE{U_{N,l}^2 \1 \{ |U_{N,l}| \geq \epsilon\}}{\mcf{N,0}} \eqsp.
\end{multline}
Using either (i) or (ii) and since $\epsilon > 0$ is arbitrary, it follows from \eqref{eq:lem:lindeberg-condition} that the \rhs\ tends in probability to $0$ as $N \to \infty$.
Using again a Taylor inequality,
\begin{multline}
\label{eq:lem:CLT:rel3}
\CPE{ \sum_{l=1}^{M_N} \left| \CPE{\exp \left( - u^2  \sigma_{N,l}^2 /2 \right) - 1 - \frac{1}{2} u^2 \sigma_{N,l}^2 }{\mcf{N,l-1}}\right|}{\mcf{N,0}} \\
\leq  \frac{u^4}{8} \sum_{l=1}^{M_N} \CPE{\sigma_{N,l}^4}{\mcf{N,0}} \leq \frac{u^4}{8} \CPE{\max_{1 \leq l \leq M_N} \sigma_{N,l}^2 }{\mcf{N,0}}  \left(\sum_{l=1}^{M_N} \sigma_{N,l}^2\right)\eqsp,
\end{multline}
since under $(i)$ or $(ii)$, $\sum_{l=1}^{M_N} \sigma_{N,l}^2$ is $\mcf{N,0}$-measurable. Because, for any $\epsilon > 0$,
\[
\CPE{\max_{1 \leq j \leq M_N} \sigma_{N,j}^2}{\mcf{N,0}} \leq \epsilon^2 + \sum_{j=1}^{M_N} \CPE{U_{N,j}^2 \1 \{ |U_{N,j}| \geq \epsilon \}}{\mcf{N,0}} \eqsp,
\]
it follows from \eqref{eq:lem:lindeberg-condition} and assumptions (i) or (ii) that the \rhs\ of \eqref{eq:lem:CLT:rel3} tends in probability to $0$ as $N \to \infty$.
Therefore the \lhs\ of \eqref{eq:lem:CLT:rel1} tends in probability to $0$ because the sum on the \rhs\ of this display is bounded by \eqref{eq:lem:CLT:rel2Bis} and \eqref{eq:lem:CLT:rel3}.
\end{proof}

\begin{prop}
\label{prop:CLTTriangularArrayLindeberg}
Assume that  $\sigma_{N,i}^2 \eqdef \CPE{U_{N,i}^2}{\mcf{N,i-1}} < \infty$ $\PP$-\as\ ,
$\CPE{U_{N,i}}{\mcf{N,i-1}}=0$ for all $N$ and $i=1, \dots, M_N$ and
\begin{align}
\label{eq:CLT:tendstoconstant}
&\sum_{i=1}^{M_N} \sigma_{N,i}^2 \plim \sigma^2 \quad && \text{for some $\sigma^2 > 0$} \eqsp, \\
\label{eq:CLT:lindeberg-condition}
&\sum_{i=1}^{M_N} \CPE{U_{N,i}^2 \1 { \{ |U_{N,i}| \geq \epsilon  \}}}{\mcf{N,i-1}} \plim 0 && \text{for any $\epsilon > 0$} \eqsp.
\end{align}
Then for any real $u$,
\begin{equation}
\label{eq:limittcl}
\CPE{\exp \left(  \rmi u \sum_{i=1}^{M_N} U_{N,i} \right)}{\mcf{N,0}} \plim \exp \left( - \sigma^2u^2/2  \right) \eqsp.
\end{equation}
\end{prop}

\begin{proof}
Without loss of generality, we assume that $\sigma^2=1$. Define the stopping time $\tau_N$  by
\[
\tau_N \eqdef \max \left\{ 1 \leq k \leq M_N: \sum_{j=1}^k \sigma_{N,j}^2 \leq 1 \right\} \eqsp,
\]
with the convention $\max \emptyset =0$. Put
$\bar{U}_{N,k}= U_{N,k}$ for $k \leq \tau_N$, $\bar{U}_{N,k}= 0$ for
$\tau_N<k \leq M_N$ and $\bar{U}_{N,M_N+1}= \left(1 -
\sum_{j=1}^{\tau_N} \sigma_{N,j}^2 \right)^{1/2} Y_N$, where $\{
Y_N\}$ are $\mathcal{N}(0,1)$ independent and independent of
$\mcf{N,M_N}$. By construction,
\[
\CPE{\bar{U}^2_{N,M_N+1}}{\mcf{N,M_N}} = 1 - \sum_{j=1}^{\tau_N} \sigma_{N,j}^2 < \infty \eqsp, \quad \PP-\as\
\]
and $\CPE{\bar{U}_{N,M_N+1}}{\mcf{N,M_N}}= 0$.
The triangular array $\{ \bar{U}_{N,k} \}_{1 \leq k \leq M_N+1}$ obviously satisfies the conditions of Lemma \ref{lem:CLT:TriangularArrayLindeberg} for
conditional means and variances. By construction,
\begin{equation}
\label{eq:prop:CLT:rel1}
\sum_{j=1}^{M_N} \CPE{\bar{U}_{N,j}^2 \1 \{ |\bar{U}_{N,j}| \geq \epsilon\}}{\mcf{N,0}} = \CPE{\sum_{j=1}^{\tau_N} \CPE{U_{N,j}^2 \1 \{ |U_{N,j}| \geq \epsilon \}}{\mcf{N,j-1}} }{\mcf{N,0}} \eqsp.
\end{equation}
Since $\sum_{j=1}^{\tau_N} \CPE{U_{N,j}^2 \1 \{ |U_{N,j}| \geq \epsilon\}}{\mcf{N,j-1}} \leq 1$, then \eqref{eq:CLT:lindeberg-condition} shows that
the \rhs\ of \eqref{eq:prop:CLT:rel1} converges in probability to $0$ as $N \to \infty$. On the other hand, if $\tau_N < M_N$,
\[
0 \leq 1 - \sum_{j=1}^{\tau_N} \sigma_{N,j}^2 \leq \sigma_{N,\tau_N+1}^2 \leq \max_{1 \leq j \leq M_N} \sigma_{N,j}^2 \eqsp.
\]
For any $\epsilon > 0$,
\[
\max_{1 \leq j \leq M_N} \sigma_{N,j}^2 \leq \epsilon^2 + \sum_{j=1}^{M_N} \CPE{U_{N,j}^2 \1 \{ |U_{N,j}| \geq \epsilon\}}{\mcf{N,j-1}}
\]
Since $\epsilon > 0$ is arbitrary, it follows from \eqref{eq:CLT:lindeberg-condition} that
\begin{equation}
\label{eq:prop:CLT:rel2}
1 - \sum_{j=1}^{\tau_N} \sigma_{N,j}^2 \plim 0 \eqsp,
\end{equation}
and since $ 1 - \sum_{j=1}^{\tau_N} \sigma_{N,j}^2 \leq 1$, this implies that $\CPE{\bar{U}^2_{N,M_N+1}}{\mcf{N,0}} \plim 0$. Therefore,
$\{ \bar{U}_{n,k} \}_{1 \leq k \leq M_N+1}$ satisfies \eqref{eq:lem:lindeberg-condition}. Put
\begin{equation}
\label{eq:prop:CLT:rel3}
\sum_{j=1}^{M_N} U_{N,j} = \sum_{j=1}^{M_N+1} \bar{U}_{N,j} - \bar{U}_{N,M_N+1}+ \sum_{j=\tau_N+1}^{M_N} U_{N,j}  \eqsp.
\end{equation}
The first term on the \rhs\ is asymptotically $\mathcal{N}(0,1)$ and $\bar{U}_{N,M_N+1} \plim 0$. It remains to prove that $\sum_{j=\tau_N+1}^{M_N} U_{N,j} \plim 0$. First note that
\begin{equation}
  \label{eq:termeRestant}
   \sum_{j=\tau_N+1}^{M_N} \sigma_{N,j}^2 =\sum_{j=1}^{M_N} \sigma_{N,j}^2-1+\left(1-\sum_{j=1}^{\tau_N} \sigma_{N,j}^2\right)\plim 0
\end{equation}
For any $\lambda > 0$,
\begin{equation*}
\PE\left( \sum_{j=\tau_N+1}^{M_N} U_{N,j} \1 \left\{ \sum_{i=\tau_N+1}^{j} \sigma_{N,i}^2 \leq \lambda \right\} \right)^2 = \PE \left(\sum_{j=\tau_N+1}^{M_N} \sigma_{N,j}^2 \1 \left\{ \sum_{i=\tau_N+1}^{j} \sigma_{N,i}^2 \leq \lambda \right\} \right) \eqsp.
\end{equation*}
The term between brackets converges to 0 in probability by (\ref{eq:termeRestant}) and its absolute value is bounded by $\lambda$. Thus, by dominated convergence, this expectation converges to 0. Thus,
$$
\sum_{j=\tau_N+1}^{M_N} U_{N,j} \1 \left\{ \sum_{i=\tau_N+1}^{j} \sigma_{N,i}^2 \leq \lambda \right\} \plim 0
$$
Moreover,
\begin{align*}
&\PP \left( \sum_{j=\tau_N+1}^{M_N} U_{N,j} \1 \left\{ \sum_{i=\tau_N+1}^j \sigma_{N,i}^2 > \lambda \right\}  \ne 0\right)\\
&\quad \leq \PP \left(\exists j \in \{1, \ldots, M_N\}, \quad \sum_{i=\tau_N+1}^j \sigma_{N,i}^2 > \lambda \right) =\PP \left(\sum_{i=\tau_N+1}^{M_N} \sigma_{N,i}^2 > \lambda \right) \eqsp,
\end{align*}
which converges to 0 by (\ref{eq:termeRestant}). The proof is completed.
\end{proof}

We will use the following technical Lemma which is a conditional version of \cite[Lemma 3.3]{dvoretzky:1972}.

\begin{lem} \label{lem:curiousLemma}
Let $\mcg{}$ a $\sigma$-field and $X$ a random variable such that $\CPE{X^2}{\mcg{}}< \infty$. Then, for any $\epsilon>0$,
  \[
4 \CPE{|X|^2 \1 \{ |X| \geq \epsilon \}}{\mcg{}} \geq \CPE{|X - \CPE{X}{\mcg{}}|^2 \1 \{ |X - \CPE{X}{\mcg{}}| \geq 2 \epsilon \} }{\mcg{}} \eqsp,
\]
\end{lem}
\begin{proof}
  Let $Y=X - \CPE{X}{\mcg{}}$. We have $\CPE{Y}{\mcg{}}=0$. It is equivalent to show that for any $\mcg{}$-measurable random variable $c$,
$$
 \CPE{Y^2 \1 \{ |Y| \geq 2 \epsilon \} }{\mcg{}} \leq 4 \CPE{|Y+c|^2 \1 \{ |Y+c| \geq \epsilon \}}{\mcg{}}
$$
On the set $\{ |c|< \epsilon\}$,
\begin{eqnarray*}
  \CPE{Y^2 \1 \{ |Y| \geq 2 \epsilon \} }{\mcg{}} & \leq &  2 \CPE{((Y+c)^2+c^2) \1 \{ |Y+c| \geq \epsilon \} }{\mcg{}}\\
& \leq &2(1+c^2/\epsilon^2) \CPE{(Y+c)^2 \1 \{ |Y+c| \geq \epsilon \} }{\mcg{}}\\
& \leq &4 \CPE{(Y+c)^2 \1 \{ |Y+c| \geq \epsilon \} }{\mcg{}}\eqsp.
\end{eqnarray*}
Moreover, on the set  $\{ |c|\geq \epsilon\}$, using that $\CPE{cY}{\mcg{}}=c\CPE{Y}{\mcg{}}=0$.
\begin{eqnarray*}
  \CPE{Y^2 \1 \{ |Y| \geq 2 \epsilon \} }{\mcg{}} & \leq &  \CPE{Y^2+c^2-\epsilon^2}{\mcg{}}\\
& \leq &\CPE{(Y+c)^2 -\epsilon^2 }{\mcg{}}\\
& \leq &\CPE{(Y+c)^2 \1 \{ |Y+c| \geq \epsilon \} }{\mcg{}}\eqsp.
\end{eqnarray*}
The proof is completed.
\end{proof}

\begin{cor}
\label{cor:CLT-NZM}
Assume that for each $N$ and $i=1, \dots, M_N$, $\CPE{U_{N,i}^2}{\mcf{N,i-1}} < \infty$ and
\begin{align}
\label{eq:CLT-NZM:tendstoconstant}
& \sum_{i=1}^{M_N} \{ \CPE{U^2_{N,i}}{\mcf{N,i-1}} - ( \CPE{U_{N,i}}{\mcf{N,i-1}})^2  \} \plim \sigma^2 && \text{for some $\sigma^2 >0$} \eqsp, \\
\label{eq:CLT-NZM:lindeberg-condition}
& \sum_{i=1}^{M_N} \CPE{U_{N,i}^2 \1_{\{ |U_{N,i}| \geq \epsilon \}}}{\mcf{N,i-1}} \plim 0  && \text{for any  $\epsilon > 0$} \eqsp.
\end{align}
Then, for any real $u$,
\begin{equation}
  \label{eq:conditionalTcl}
  \CPE{\exp \left( \rmi u \sum_{i=1}^{M_N} \{ U_{N,i} - \CPE{U_{N,i}}{\mcf{N,i-1}} \} \right)}{\mcf{N,0}} \plim \exp( - (u^2/2) \sigma^2) \eqsp.
\end{equation}
If in addition $\mcf{N,i} \subseteq \mcf{N',i}$ for all $N'\geq N$ and $i \leq M_N$, then
\begin{equation}
  \label{eq:mixing}
  \sum_{i=1}^{M_N} \{ U_{N,i} - \CPE{U_{N,i}}{\mcf{N,i-1}} \}  \dlim {\mathcal N}(0, \sigma^2) \quad (mixing)\eqsp.
\end{equation}
\end{cor}

\begin{proof}
Set $\bar{U}_{N,i} = U_{N,i} - \CPE{U_{N,i}}{\mcf{N,i-1}}$. By construction
$\CPE{\bar{U}_{N,i}}{\mcf{N,i-1}}= 0$ and because $\CPE{\bar{U}_{N,i}^2}{\mcf{N,i-1}} =\CPE{U_{N,i}^2}{\mcf{N,i-1}} - ( \CPE{U_{N,i}}{\mcf{N,i-1}} )^2$, \eqref{eq:CLT:tendstoconstant} is fulfilled. The proof of \eqref{eq:CLT:lindeberg-condition} follows from Lemma \ref{lem:curiousLemma}.
It remains to show (\ref{eq:mixing}). Let $\mcf{\infty}=\vee_{1}^\infty \mcf{N,M_N}$ the $\sigma$-field generated by $ \cup_1^\infty \mcf{N,M_N}$. For any $E' \in \mcf{\infty}$ and any $ \epsilon >0$, there exists an $m$ such that $E \in \mcf{m,M_m}$ and $\PP(E \Delta E')< \epsilon$ (where $\Delta$ denotes the symmetric difference). Denotes $V_{N,i}= U_{N,i} - \CPE{U_{N,i}}{\mcf{N,i-1}} $.
$$
\left|\PE\left[\exp\left(\rmi u \sum_{i=1}^{M_N} V_{N,i} \right) \1(E)\right]-\PE\left[\exp\left(\rmi u \sum_{i=1}^{M_N} V_{N,i}  \right)\1(E')\right]\right| \leq \PP(E \Delta E')< \epsilon
$$
Moreover, applying (\ref{eq:conditionalTcl}) to the triangular array $\{U_{N,M_m+i}\}_{1 \leq i \leq M_N -M_m}$ (with $N$ sufficiently large) associated to the $\sigma$-field $\{\mcf{N,M_m+i}\}_{0 \leq i \leq M_N -M_m}$,
$$
\PE\left[\exp\left(\rmi u \sum_{i=m+1}^{M_N} V_{N,i}  \right)\1(E)\right] \plim \exp(-u^2\sigma^2/2)\PP(E)
$$
and $\max_{1 \leq i \leq M_N} |V_{N,i}| \plim 0$. Thus,
$$
\PE\left[\exp\left(\rmi u \sum_{i=1}^{M_N} V_{N,i}  \right)\1(E)\right] \plim \exp(-u^2\sigma^2/2)\PP(E)
$$
This implies that for any  $E' \in \mcf{\infty}$,
$$
\PE\left[\exp\left(\rmi u \sum_{i=1}^{M_N} V_{N,i}  \right)\1(E')\right] \plim \exp(-u^2\sigma^2/2)\PP(E')
$$
Now, let $E \in \mcf{}$,
\begin{multline*}
  E\left[\exp\left(\rmi u \sum_{i=1}^{M_N} V_{N,i}  \right)\1(E)\right]= \PE\left[\exp\left(\rmi u \sum_{i=1}^{M_N} V_{N,i}  \right)\CPE{\1(E)}{\mcf{\infty}}\right]\\
 \plim  \exp(-u^2\sigma^2/2)\PE(\CPE{\1(E)}{\mcf{\infty}})=  \exp(-u^2\sigma^2/2)\PP(E).
\end{multline*}
The proof is completed.
\end{proof}

\section{Proof of  Theorems \ref{thm:LLN-mutation} and \ref{thm:CLT-mutation}}
\begin{proof}[Proof of the Theorem \ref{thm:LLN-mutation}]
We set $\mcf{N,0} = \sigma \left( \left\{ (\epart{N,i}, \ewght{N,i}) \right\}_{1 \leq i \leq M_N}  \right)$ and for $j=1, \dots, \tilde{M}_N$,
$\mcf{N,j}= \mcf{N,0} \vee \sigma \left( \{ \etpart{N,k} \}_{1 \leq k \leq j} \right)$.
Checking that $\tilde{\consistfunc}$ is proper is straightforward, so
we turn to the consistency. We show first that for any $f \in \tilde{\consistfunc}$,
\begin{equation}
\frac{1}{\alpha_N \sumweight{N}} \sum_{j=1}^{\tilde{M}_N} \etwght{N,j} f(\etpart{N,j}) \plim
\nu  \LTAR( f) \eqsp, \label{eq:limit-limit}
\end{equation}
where $\etpart{N,j}$ and $\etwght{N,j}$  are defined in \eqref{eq:mutationstep-1} and \eqref{eq:mutationstep-3}, respectively.
Because the mutation step is unbiased (see \eqref{eq:unbiased-mutation}),
\[
(\alpha_N \sumweight{N})^{-1}
\sum_{j=1}^{\tilde{M}_N} \CPE{\etwght{N,j} f(\etpart{N,j})}{\mcf{N,j-1}}
=\sumweight{N}^{-1} \sum_{i=1}^{M_N} \ewght{N,i} \LTAR (\epart{N,i},f) \eqsp.
\]
Because the weighted sample $\{ (\epart{N,i}, \ewght{N,i}) \}_{1 \leq i \leq M_N}$ is consistent for $(\nu,\consistfunc)$ and for
$f \in \tilde{\consistfunc}$, the function $\LTAR(\cdot,f)  \in \consistfunc$,
\[
\sumweight{N}^{-1} \sum_{i=1}^{M_N} \ewght{N,i} \LTAR (\epart{N,i},f) \plim  \nu \LTAR (f)\eqsp,
\]
it suffices to show that
\begin{equation}
\label{eq:consistency-weighted:keyresult}
(\alpha_N \sumweight{N})^{-1}\sum_{j=1}^{\tilde{M}_N} \left \{  \etwght{N,j} f(\etpart{N,j})
-  \CPE{\etwght{N,j} f(\etpart{N,j})}{\mcf{N,j-1}} \right \} \plim 0 \eqsp.
\end{equation}
Put $U_{N,j}= (\alpha_N \sumweight{N})^{-1}\etwght{N,j} f(\etpart{N,j})$ for $j=1, \dots, \tilde{M}_N$
and appeal to Proposition~\ref{prop:LLN-NZM}. Just as above,
\[
\sum_{j=1}^{\tilde{M}_N} \CPE{|U_{N,j}|}{\mcf{N,j-1}}
= \sumweight{N}^{-1}\sum_{i=1}^{M_N} \ewght{N,i} \LTAR(\epart{N,i},|f|)
\plim \nu \LTAR(|f|) \eqsp,
\]
showing that the sequence $ \left \{ \sum_{j=1}^{\tilde{M}_N} \CPE{|U_{N,j}|}{\mcf{N,j-1}} \right \}_{N \geq 0}$
is tight (Proposition~\ref{prop:LLN-NZM}-Eq.\eqref{eq:LLN-NZM:tightness}).
For any $\epsilon > 0$, put $A_N \eqdef \sum_{j=1}^{\tilde{M}_N} \CPE{| U_{N,j} | \1 {\{ |U_{N,j}| \geq \epsilon \}}}{\mcf{N,j-1}}$.
We  need to show that $ A_N  \plim 0$ (Proposition~\ref{prop:LLN-NZM}-Eq.\eqref{eq:LLN-NZM:negligeability}).
For any positive $C$, $\epart{} \in \Xi$,
$\KISS \left(\epart{}, W |f| \1 { \{ W |f| \geq C \}}\right) \leq \KISS(\epart{},W |f|) =  \LTAR(\epart{},|f|)$.
Because the function  $\LTAR\left(\cdot,|f|\right)$ belongs to the proper set $\consistfunc$,
the function $ \KISS(\cdot, W |f| \1{ \{W |f| \geq C \}})$ belongs  to $\consistfunc$. Hence for all $C,\ \epsilon>0$,
\begin{align*}
&A_N \1 \left\{ (\alpha_N\sumweight{N})^{-1} \max_{1 \leq i \leq M_N} \ewght{N,i}  \leq  \epsilon /C \right\}\\
&\quad \leq \sumweight{N}^{-1} \sum_{i=1}^{M_N}  \ewght{N,i} \left[ \alpha_N^{-1} \sum_{k=1}^{\alpha_N} \KISS\left(\epart{N,i},  W |f| \1{\{W |f| \geq C\}}\right) \right] \\
&\quad \plim \nu \KISS \left(W |f|  \1{\{ W |f| \geq C\}} \right) \eqsp.
\end{align*}
By dominated convergence, the \rhs\ can be made arbitrarily small by letting $C \to \infty$.
Combining with $ \sumweight{N}^{-1} \max_{1 \leq i \leq M_N} \ewght{N,i}  \plim 0$, this shows that $A_N$ tends to zero in probability,
showing \eqref{eq:LLN-NZM:negligeability}.
Thus Proposition~\ref{prop:LLN-NZM} applies and \eqref{eq:limit-limit} holds.
Under the stated assumptions, the function $\LTAR(\cdot,\etpartset)$
belongs to $\consistfunc$, implying that the constant function $ g \equiv 1$
satisfies \eqref{eq:limit-limit}; therefore,
$(\alpha_N \sumweight{N})^{-1} \sum_{j=1}^{\tilde{M}_N} \etwght{N,j} \plim \nu \LTAR(\epartset)$. Combined with
\eqref{eq:limit-limit} this shows that for any $f \in \tilde{\consistfunc}$,
\[
\tsumweight{N}^{-1} \sum_{j=1}^{\tilde{M}_N} \etwght{N,j} f(\etpart{N,j}) \plim \mu(f) \eqsp,
\]
To complete the proof of the consistency, it remains to prove that
\begin{equation*}
\tsumweight{N}^{-1} \, \max_{1 \leq j \leq \tilde{M}_N}  \etwght{N,j} \plim 0 \eqsp.
\end{equation*}
Because $(\alpha_N \sumweight{N})^{-1} \tsumweight{N} \plim \nu \LTAR(\epartset)$, it is actually sufficient to show that
$$ (\alpha_N \sumweight{N})^{-1} \, \max_j \etwght{N,j} \plim 0 \eqsp.$$
For any $C>0$,
\begin{align*}
\left(\alpha_N \sumweight{N} \right)^{-1} \max_{ 1 \leq j \leq \tilde{M}_N} \etwght{N,j}\1_{\{W(\etpart{N,j})\leq C\}} &\leq C \left( \alpha_N\sumweight{N} \right)^{-1} \max_{1 \leq i \leq M_N} \ewght{N,i}  \plim 0\\
\left( \alpha_N \sumweight{N} \right)^{-1} \max_{1 \leq j \leq \tilde{M}_N} \etwght{N,j}\1_{\{W(\etpart{N,j})> C\}}  &\leq   \left(\alpha_N \sumweight{N} \right)^{-1}
\sum_{j=1}^{\tilde{M}_N} \etwght{N,j}\1_{\{W(\etpart{N,j})> C\}} \\
&\plim \nu \LTAR \left(\{W> C\}\right)
\end{align*}
The term in the RHS of the last equation goes to zero as $C \to \infty$  which concludes the proof.
\end{proof}

\begin{proof}[Proof of the Theorem \ref{thm:CLT-mutation}]
First we note that by definition $\alpha$ is necessarily at least 1.
Checking that $\tilde{\asnormfunc}$ and $\tilde{\consistweight}$ is proper is straightforward.
Pick $f \in \tilde{\asnormfunc}$ and assume, without loss of generality, that $\mu(f)= 0$.
Write
$
\tsumweight{N}^{-1} \sum_{i=1}^{\tilde{M}_N} \etwght{N,i} f(\etpart{N,i})  = (\alpha_N \sumweight{N} / \tsumweight{N})  \left(A_N +  B_N \right) \eqsp,
$
with
\begin{align*}
&A_N = (\alpha_N \sumweight{N})^{-1}\sum_{j=1}^{\tilde{M}_N} \CPE{\etwght{N,j}  f(\etpart{N,j})}{\mcf{N,j-1}}   = \sumweight{N}^{-1}\sum_{i=1}^{M_N} \ewght{N,i} \LTAR(\epart{N,i},f)   \eqsp, \\
&B_N = (\alpha_N \sumweight{N})^{-1} \sum_{j=1}^{\tilde{M}_N} \left\{ \etwght{N,j}  f(\etpart{N,j})  - \CPE{\etwght{N,j} f(\etpart{N,j})}{\mcf{N,j-1}} \right\} \eqsp.
\end{align*}
Because $\alpha_N \sumweight{N}  / \tsumweight{N}$
converges to $1 / \nu\LTAR(\etpartset)$ in probability
(see in the proof of Theorem \ref{thm:LLN-mutation}),
the conclusion of the theorem follows from
Slutsky's theorem if we prove that
$M_N^{1/2} (A_N + B_N)$ converges weakly to
$\gauss(0,\sigma^2(\LTAR f ) + \alpha^{-1} \eta^2(f))$ where
\begin{equation}
\label{eq:definition:eta2}
\eta^2(f) \eqdef \gamma  \KISS \{ [Wf - \KISS(\cdot,Wf)]^2 \} \eqsp,
\end{equation}
with $W$ given in \eqref{eq:definitionW}. The function $\LTAR(\cdot,f)$ belongs to $\asnormfunc$ and
$ \nu \LTAR(f)=  \mu(f) \, \nu\LTAR(\etpartset)=0$.
Because $\{ (\epart{N,i},\ewght{N,i}) \}_{1 \leq i \leq M_N}$
is asymptotically normal for $(\nu,\asnormfunc,\consistweight, \stdnormfunc, \gamma, \{M^{1/2}_N\})$,
$ M^{1/2}_N A_N \dlim \gauss ( 0, \stdnormfunc^2 ( \LTAR f ) )$. Next we prove that for any real $u$,
$$
\CPE{ \exp( \rmi u \tilde{M}_N^{1/2} B_N)}{\mcf{N,0}} \plim \exp \left( - (u^2/2) \eta^2(f) \right) \eqsp.
$$
where $\eta^2(f)$ is defined in \eqref{eq:definition:eta2}.
For that purpose we use corollary~\ref{cor:CLT-NZM}, and
we thus need to check \eqref{eq:CLT-NZM:tendstoconstant}-\eqref{eq:CLT-NZM:lindeberg-condition}
with
\[
U_{N,j} \eqdef (\alpha_N \sumweight{N})^{-1} \, \tilde{M}_N^{1/2} \etwght{N,j} f(\etpart{N,j})  \eqsp, \quad j=1, \dots, \tilde{M}_N \eqsp.
\]
Under the stated assumptions, for $f \in \asnormfunc$, the function $\KISS(\cdot,W^2 f^2) $
belongs to $\consistweight$. By the Jensen inequality, for any $\epart{} \in \epartset$,
$\left\{ \LTAR(\epart{},f) \right\}^2 = \left \{ \KISS(\epart{}, Wf) \right\}^2 \leq \KISS(\epart{}, W^2 f^2)$.
Because the $\consistweight$ is proper and the function $\KISS(\cdot, W^2f^2) \in \consistweight$, the relation
$\{ \LTAR(\cdot,f) \}^2 \leq \KISS(\cdot, W^2f^2)$
implies that the function $\{ \LTAR(\cdot, f )\}^2$ also belongs to $\consistweight$.
Because $\{ (\epart{N,i}, \ewght{N,i}) \}_{1 \leq i \leq M_N}$ is asymptotically normal for
$(\nu,$ $\asnormfunc,$ $\consistweight,$ $\stdnormfunc,$ $\gamma,$ $\{M^{1/2}_N\})$, \eqref{eq:AN2} implies
\[
  \sum_{j=1}^{\tilde{M}_N} \CPE{U_{N,j}^2}{\mcf{N,j-1}}= \frac{M_N}{ \sumweight{N}^{2}}\sum_{i=1}^{M_N} \ewght{N,i}^2 \KISS(\epart{N,i},W^2 f^2)
 \plim \gamma \,\KISS(W^2 f^2) \eqsp,
\]
\[
\sum_{j=1}^{\tilde{M}_N} ( \CPE{U_{N,j}}{\mcf{N,j-1}} )^2= \frac{M_N}{ \sumweight{N}^{2}}\sum_{i=1}^{M_N} \ewght{N,i}^2 \left \{\LTAR (\epart{N,i},f) \right \}^2\plim
\int \gamma(d \epart{}) \left\{ \LTAR (\epart{},f) \right\}^2\eqsp.
\]
These displays imply that \eqref{eq:CLT-NZM:tendstoconstant}
holds. It remains to check \eqref{eq:CLT-NZM:lindeberg-condition}.
We have, for all $C,\ \epsilon >0$,
\begin{multline*}
\left(\sum_{j=1}^{\tilde{M}_N} \CPE{U_{N,j}^2 \1_{\{ |U_{N,j}|\geq \epsilon \} }}{\mcf{N,j-1}}\right) \1{\left\{\frac{M_N \max_{1 \leq i \leq M_N} \ewght{N,i}^2}{\alpha_N \sumweight{N}^2}\leq \left(\frac{\epsilon}{C}\right)^2\right\}} \leq \\
\frac{M_N}{ \sumweight{N}^{2}}\sum_{i=1}^{M_N} \ewght{N,i}^2 \KISS\left(\epart{N,i}, W^2 f^2 \1 {\{|W f| \geq C\}}\right)
\plim \gamma \KISS\left(W^2 f^2 \1 {\{ |W f| \geq C \}}\right)
\end{multline*}
which converges to 0 as $C$ goes to infinity. Combining with $$ (M_N)^{1/2} \sumweight{N}^{-1} \max_{1 \leq i \leq M_N} \ewght{N,i} \plim 0$$ yields
\[
  \sum_{j=1}^{\tilde{M}_N} \CPE{U_{N,j}^2 \1 {\{ |U_{N,j}| \geq \epsilon \} }}{\mcf{N,j-1}}\plim 0 \eqsp ,
\]
and this is condition \eqref{eq:CLT-NZM:lindeberg-condition}. Thus corollary \ref{cor:CLT-NZM} applies. It follows that
$$
\left(
\begin{array}{c}
M^{1/2}_N A_N \\
\tilde{M}_N^{1/2} B_N
\end{array}
\right)
\dlim \gauss\left( 0, \mathrm{diag}\{ \sigma^2(\LTAR f ), \eta^2(f)\} \right)\eqsp,
$$
where $\eta^2(f)$ is given in \eqref{eq:definition:eta2}.
The proof of condition (\ref{eq:AN1}) in the definition of an asymptotically normal sample is now concluded upon writing
$M^{1/2}_N (A_N + B_N)= M^{1/2}_N A_N + \alpha^{-1/2}_N \tilde{M}_N^{1/2} B_N$.

Consider now \eqref{eq:AN2}. Recalling that $\tsumweight{N}/(\alpha_N \sumweight{N}) \plim \nu \LTAR(\epartset)$, it is sufficient to show that
for $f \in \tilde{\consistweight}$,
\begin{equation}
\label{eq:lgnSquareWght}
\frac{M_N}{(\alpha_N \sumweight{N})^2} \sum_{j=1}^{\tilde{M}_N} \etwght{N,j}^2 f(\etpart{N,j}) \plim
\alpha^{-1}\gamma \KISS(h) \eqsp,
\end{equation}
where $h \eqdef W^2 f$. Define $U_{N,j}= (\alpha_N \sumweight{N})^{-2} M_N \etwght{N,j}^2 f(\etpart{N,j}) $.
Under the stated assumptions, for any $f \in \tilde{\consistweight}$,
the functions $\KISS(\cdot,|h|)$ and $\KISS(\cdot,h)$ belong to $\consistweight$. Because $\{ (\epart{N,i}, \ewght{N,i}) \}_{1 \leq i \leq M_N}$
is asymptotically normal for $(\nu,\asnormfunc,\consistweight,\stdnormfunc,\gamma,\{M_N^{1/2}\})$,
\begin{align}
\label{eq:rel1}
&\sum_{j=1}^{\tilde{M}_N} \CPE{|U_{N,j}|}{\mcf{N,j-1}}
=  \frac{M_N}{\alpha_N \sumweight{N}^2} \sum_{i=1}^{M_N} \ewght{N,i}^2 \KISS(\epart{N,i},|h|)
\plim \alpha^{-1}\gamma \KISS(|h|) \eqsp, \\
\label{eq:rel2}
&\sum_{j=1}^{\tilde{M}_N} \CPE{U_{N,j}}{\mcf{N,j-1}}
= \frac{M_N}{\alpha_N \sumweight{N}^2} \sum_{i=1}^{M_N} \ewght{N,i}^2 \KISS(\epart{N,i},h) \plim  \alpha^{-1}\gamma \KISS(h)\eqsp.
\end{align}
We appeal to Proposition \ref{prop:LLN-NZM}. Eq.~\eqref{eq:rel1} shows that the tightness condition \eqref{eq:LLN-NZM:tightness}. For $\epsilon > 0$, set
$A_N= \sum_{j=1}^{\tilde{M}_N} \CPE{|U_{N,j}| \1_{\{|U_{N,j}| \geq \epsilon\}}}{\mcf{N,j-1}}$.
We  need to show that $A_N \plim 0$. For any positive $C$ and $\epart{} \in \epartset$,
$ \KISS\left(\epart{},|h| \1 { \{ |h| \geq C \}}\right) \leq \KISS(\epart{},|h|)$. Because the function $\KISS(\cdot,|h|)$ belongs to $\consistweight$ and
the set $\consistweight$ is proper, the function $\KISS(\cdot, |h| \1{ \{ |h| \geq C \}})$ belongs to $\consistfunc$. Hence for all $C,\ \epsilon>0$,
\begin{align*}
&A_N \1 {\left\{\max_i\frac{M_N \ewght{N,i}^2}{ (\alpha_N \sumweight{N})^2}\leq \frac{\epsilon}{C}\right\}} \\
\quad &\leq \frac{M_N}{ \sumweight{N}^2} \sum_{i=1}^{M_N}  \ewght{N,i}^2 \KISS\left(\epart{N,i}, | h| \1 {\{|h| \geq C\}}\right)
\plim \gamma \KISS\left(|h| \1 {\{ |h| \geq C\}}\right) \eqsp.
\end{align*}
By the dominated convergence theorem, the \rhs\ can be made arbitrarily small by letting $C \to \infty$.
Combining with $ M_N (\alpha_N \sumweight{N}^2)^{-1} \, \max_{ 1 \leq i \leq M_N} \ewght{N,i}^2   \plim 0$, this shows that $A_N$ tends to zero in probability,
showing \eqref{eq:LLN-NZM:negligeability}.
Thus Proposition~\ref{prop:LLN-NZM} applies and condition (\ref{eq:AN2}) is proved. To complete the proof, it remains to prove (\ref{eq:AN3}).
Combining with $\tsumweight{N}/(\alpha_N \sumweight{N})\plim \nu \LTAR(\epartset)$ (see proof of Theorem \ref{thm:LLN-mutation})
and $\tilde{M}_N=\alpha_N M_N$, it is actually sufficient to show that
$ \left( \alpha_N \sumweight{N}^2 \right)^{-1} M_N \max_{1 \leq j \leq \tilde{M}_N} \etwght{N,j}^2  \plim 0$.
For any $C>0$,
\begin{align*}
M_N  \frac{\max_{1 \leq j \leq \tilde{M}_N}
  \etwght{N,j}^2\1_{\{W(\etpart{N,j})\leq C\}}}{\alpha_N (\sumweight{N})^2} \leq C^2  M_N \frac{\max_{1 \leq i \leq M_N}
  \ewght{N,i}^2}{\alpha_N(\sumweight{N})^2} \plim 0 \eqsp.
\end{align*}
Applying \eqref{eq:lgnSquareWght} with $f \equiv 1$,  it holds that
\begin{multline*}
M_N \frac{\max_{1 \leq j \leq \tilde{M}_N} \etwght{N,j}^2\1_{\{W(\etpart{N,j})> C\}}}{ \alpha_N (\sumweight{N})^2} \leq
\frac{ M_N}{\alpha_N (\sumweight{N})^2} \sum_{j=1}^{\tilde{M}_N} \etwght{N,j}^2\1_{\{W(\etpart{N,j}) \geq C\}}\\
 \plim \gamma \KISS \left(W^2  \1{\{W \geq C\}} \right) \eqsp.
\end{multline*}
The term on the \rhs\ goes to zero as $C \to \infty$. This proves condition (\ref{eq:AN3}) and concludes the proof.
\end{proof}

\section{Proof of the Theorems \ref{thm:LLN-unbiased-sampling} and \ref{thm:CLT-independent-sampling}}

\begin{proof}[Proof of the Theorem \ref{thm:LLN-unbiased-sampling}]
As above, we set $\mcf{N,0} = \sigma \left( \left\{ (\epart{N,i}, \ewght{N,i}) \right\}_{1 \leq i \leq M_N}  \right)$ and for $j=1, \dots, \tilde{M}_N$,
$\mcf{N,j}= \mcf{N,0} \vee \sigma \left( \{ \etpart{N,k} \}_{1 \leq k \leq j} \right)$.
Pick $f$ in $\consistfunc$.
Since $\consistfunc$ is proper,
$|f| \1 { \{ |f| \geq C \}} \in \consistfunc$ for any $C \geq 0$.
Because $\{ (\epart{N,i},$ $\ewght{N,i}) \}_{1 \leq i \leq M_N}$ is consistent for $(\nu,\consistfunc)$,
and $\consistfunc$ is a proper set of functions,
\begin{multline} \label{eq:limUnbiaised}
\tilde{M}_N^{-1} \sum_{i=1}^{\tilde{M}_N}
\CPE{|f(\etpart{N,i})| \1_{\{ |f(\etpart{N,i})| \geq C \}}}{\mcf{N,i-1}}
\\
= \sumweight{N}^{-1} \sum_{i=1}^{M_N} \ewght{N,i}  |f(\epart{N,i})| \1_{\{|f(\epart{N,i})| \geq C \}} \plim \nu( |f| \1_{\{|f| \geq C\}}) \eqsp,
\end{multline}
We now check \eqref{eq:LLN-NZM:tightness}--\eqref{eq:LLN-NZM:negligeability} of  Proposition~\ref{prop:LLN-NZM}.
For any $i=1, \dots, M_N$, $U_{N,i} \eqdef \tilde{M}_N^{-1} f(\etpart{N,i})$ . Taking $C= 0$ in (\ref{eq:limUnbiaised}):
$$
\sum_{i=1}^{\tilde{M}_N} \CPE{|U_{N,i}|}{\mcf{N,i-1}}= \tilde{M}_N^{-1} \sum_{i=1}^{\tilde{M}_N} \CPE{|f(\etpart{N,i})|}{\mcf{N,i-1}} \plim \nu(|f|) < \infty \eqsp,
$$
whence the sequence
$\{ \sum_{i=1}^{\tilde{M}_N} \CPE{|U_{N,i}|}{\mcf{N,i-1}} \}_{N \geq 0}$ is tight.
Next, for any positive $\epsilon$ and $C$ we have for sufficiently large $N$,
\begin{multline*}
\sum_{i=1}^{\tilde{M}_N} \CPE{|U_{N,i}| \1_{\{ |U_{N,i}| \geq \epsilon \}}}{\mcf{N,i-1}}
=  \frac{1}{\tilde{M}_N} \sum_{i=1}^{\tilde{M}_N} \CPE{|f(\etpart{N,i})| \1_{ \{ |f|(\etpart{N,i}) \geq \epsilon \tilde{M}_N \}}}{\mcf{N,i-1}} \\
\leq  \tilde{M}_N^{-1} \sum_{i=1}^{\tilde{M}_N} \CPE{|f(\etpart{N,i})| \1_{ \{ |f|(\etpart{N,i}) \geq C \}}}{\mcf{N,i-1}} \plim \mu \left( |f|\1_{ \{ |f| \geq C \}} \right) \eqsp.
\end{multline*}
By dominated convergence the \rhs\ of this display tends to zero as
$C\to\infty$. Thus the \lhs\ of the display converges to zero in
probability, showing \eqref{eq:LLN-NZM:negligeability}.
\end{proof}

\begin{proof}[Proof of the Theorem \ref{thm:CLT-independent-sampling}]
Pick $f \in \tilde{\asnormfunc}$ and write
$\tilde{M}_N^{-1} \sum_{i=1}^{\tilde{M}_N} f(\etpart{N,i}) - \nu(f)
=A_N+B_N$ with
\begin{align*}
A_N &= \sumweight{N}^{-1} \sum_{i=1}^{M_N} \ewght{N,i}  \{ f(\epart{N,i}) - \nu(f) \} \eqsp, \\
B_N &= \tilde{M}_N^{-1} \sum_{i=1}^{\tilde{M}_N} \left \{ f(\etpart{N,i}) -  \CPE{f(\etpart{N,i})}{\mcf{N,i-1}} \right \} \eqsp,
\end{align*}
We first prove that
\begin{equation}
\CPE{ \exp\left( \rmi u \; \tilde{M}_N^{1/2} \; B_N  \right) }{\mcf{N,0}}
\plim \; \exp\left(- (u^2/2)\PVar_{\nu}(f)\right) \eqsp. \label{eq:cltcond}
\end{equation}
We will appeal to corollary \ref{cor:CLT-NZM} and hence need
to check \eqref{eq:CLT-NZM:tendstoconstant}-\eqref{eq:CLT-NZM:lindeberg-condition} with $U_{N,i} \eqdef \tilde{M}_N^{-1/2} f(\etpart{N,i})$.
First, because $\{ (\epart{N,i},\ewght{N,i}) \}_{1 \leq i \leq M_N}$ is consistent for $(\nu,\consistfunc)$ and for $f \in \tilde{\asnormfunc}$,
$f^2 \in \consistfunc$,
\begin{multline*}
\sum_{j=1}^{\tilde{M}_N} \{ \CPE{U^2_{N,j}}{\mcf{N,j-1}} - ( \CPE{U_{N,j}}{\mcf{N,j-1}})^2  \} \\=  \sumweight{N}^{-1} \sum_{i=1}^{M_N} \ewght{N,i}  f^2(\epart{N,i}) - \left \{ \sumweight{N}^{-1} \sum_{i=1}^{M_N} \ewght{N,i} f(\epart{N,i}) \right \}^2
\plim \nu(f^2) - \{ \nu(f) \}^2 = \PVar_\nu(f) \eqsp,
\end{multline*}
showing \eqref{eq:CLT-NZM:tendstoconstant}.
Pick $\epsilon > 0$. For any positive constant $C$,
\begin{multline*}
\sum_{j=1}^{\tilde{M}_N} \CPE{U_{N,j}^2 \1_{\{ |U_{N,j}| \geq \epsilon \}}}{\mcf{N}} \leq
\tilde{M}_N^{-1} \sum_{i=1}^{\tilde{M}_N} \CPE{f^2(\etpart{N,j}) \1 {\{ |f(\etpart{N,j})| \geq C \}}}{\mcf{N}} \\
= \sumweight{N}^{-1} \sum_{i=1}^{M_N} \ewght{N,i} f^2(\epart{N,i}) \1 { \{ |f|(\epart{N,i}) \geq C \}} \eqsp,
\end{multline*}
where the inequality holds for sufficiently large $N$.
Since $f^2 $ belongs to the proper set $\consistfunc \subseteq \lone(\epartset,\nu)$,
we have  $f^2 \1 { \{ |f| \geq C \}}\in \consistfunc$.
This implies that the \rhs\ of the above display converges in
probability to $\nu( f^2 \1 {\{ |f| \geq C \}})$. Because $f^2 \in \consistfunc \subseteq \lone(\epartset,\nu)$, $\nu( f^2 \1 {\{|f| \geq C\}})$
tends to zero as $C \to \infty$, so that \eqref{eq:CLT-NZM:lindeberg-condition} is satisfied.

Combining (\ref{eq:cltcond}) with $a_N A_N \dlim \gauss \left( 0, \sigma^2(f) \right)$,  we find that for any
real numbers $u$ and $v$,
\begin{multline*}
\PE\left[ \exp( \rmi (u \tilde{M}^{1/2}_N  B_N + v a_N A_N)\right] = \PE \left[
\CPE{\exp( \rmi u \tilde{M}^{1/2}_N B_N )}{\mcf{N,0}} \exp( \rmi v a_N A_N) \right] \\
\to \exp \left\{ - (u^2/2) \PVar_\nu(f) \right\} \exp \left\{ - (v^2/2) \sigma^2 \left(f \} \right) \right\} \eqsp.
\end{multline*}
Thus the bivariate characteristic function converges to the characteristic
function of a bivariate normal, implying that
$$
\left(
\begin{array}{c}
a_N A_N \\
\tilde{M}_N^{1/2} B_N
\end{array}
\right)
\dlim
\gauss\left( 0, \mathrm{diag} \left[ \sigma^2 \left( f \right), \PVar_\mu(f) \right] \right)\eqsp.
$$
Put $b_N = a_N$ if $\alpha < 1$ and $b_N = \tilde{M}_N^{1/2}$ if $\alpha \geq 1$. The proof follows from
$$
b_N (A_N + B_N) = (b_N a^{-1}_N) a_N A_N + (b_N \tilde{M}_N^{-1/2}) \tilde{M}^{1/2}_N B_N.
$$
The condition (\ref{eq:AN2}) and (\ref{eq:AN3}) are obviously fulfilled using that the weighted sample  $\{ (\etpart{N,i},1) \}_{1 \leq i \leq \tilde{M}_N}$
is consistent for $(\nu, \consistfunc)$.
\end{proof}
\section{Proof of Theorem \ref{thm:resid}}
\label{sec:proof:thm:resid}
\begin{proof}[Proof of Theorem \ref{thm:resid}]
To apply Corollary \ref{cor:CLT-NZM}, we just have to check (\ref{eq:CLT-NZM:tendstoconstant}) and (\ref{eq:CLT-NZM:lindeberg-condition}) where
$U_{N,i}= a_N \tilde{M}_N^{-1} f(\etpart{N,i})$ and $\{\mcf{N,i}\}$ defined by $\mcf{N,0}=\sigma\{(\epart{N,i})_{1 \leq i \leq M_N}\}$
and for all $1\leq k \leq \tilde{M}_N$, $\mcf{N,k}=\mcf{N,0} \vee \sigma\{(\etpart{N,i})_{1 \leq i \leq k}\}$.
Noting that $\etpart{N,j}$ is $\mcf{N,0}$ measurable for $j=1, \dots, \bar{M}_N$, we have
\begin{eqnarray}
A_N&=&  \sum_{i=1}^{\tilde{M}_N} \{ \CPE{U^2_{N,i}}{\mcf{N,i-1}} - ( \CPE{U_{N,i}}{\mcf{N,i-1}})^2  \}\nonumber \\
&= &\sum_{i=\bar{M}_N+1}^{\tilde{M}_N} \{ \CPE{U^2_{N,i}}{\mcf{N,i-1}} - ( \CPE{U_{N,i}}{\mcf{N,i-1}})^2  \} \nonumber\\
&= &\frac{a_N^2}{\tilde{M}_N} \frac{\tilde{M}_N-\bar{M}_N}{\tilde{M}_N} \left\{\sum_{i=1}^{M_N} \etwght{N,i} f^2(\epart{N,i})-\left(\sum_{i=1}^{M_N} \etwght{N,i} f(\epart{N,i})\right)^2\right\} \label{eq:condvarRes}
\end{eqnarray}
where the weights $\etwght{N,i}$ are given in \eqref{eq:multinomial-sampling-resid-2}. Note that
\begin{align*}
&\etwght{N,i}= \frac{\sumweight{N}^{-1} \ewght{N,i} -\tilde{M}_N^{-1} \lfloor \tilde{M}_N  \sumweight{N}^{-1} \ewght{N,i} \rfloor}
{1- \tilde{M}_N^{-1} \sum_{i=1}^{M_N} \lfloor \tilde{M}_N \sumweight{N}^{-1} \ewght{N,i}\rfloor } \quad \mbox{and} \quad
\frac{\tilde{M}_N-\bar{M}_N}{\tilde{M}_N}= 1-  \frac{\sum_{i=1}^{M_N} \lfloor \tilde{M}_N \sumweight{N}^{-1} \ewght{N,i}  \rfloor}{\tilde{M}_N} \eqsp.
\end{align*}
By applying  Lemma \ref{lem:convRes},
\begin{multline*}
A_N = \frac{a_N^2}{\tilde{M}_N} \sum_{i=1}^{M_N} \left(\sumweight{N}^{-1} \ewght{N,i} -\tilde{M}_N^{-1} \lfloor \tilde{M}_N  \sumweight{N}^{-1} \ewght{N,i} \rfloor \right) f^2(\epart{N,i})- \\
\frac{a_N^2}{\tilde{M}_N} \frac{\left(\sum_{i=1}^{M_N} \left( \sumweight{N}^{-1} \ewght{N,i} -\tilde{M}_N^{-1} \lfloor \tilde{M}_N  \sumweight{N}^{-1} \ewght{N,i} \rfloor \right) f(\epart{N,i})\right)^2}{1- \tilde{M}_N^{-1} \,\sum_{i=1}^{M_N} \lfloor \tilde{M}_N \sumweight{N}^{-1} \ewght{N,i}  \rfloor }
\plim  \tilde \stdnormfunc^2(f),
\end{multline*}
It remains to check \eqref{eq:AN2} and \eqref{eq:AN3}.
By Theorem \ref{thm:LLN-unbiased-sampling}, the weighted sample $\{(\etpart{N,i},1)\}_{ } $ is consistent for $(\nu,\consistfunc)$, which implies,
$$
\frac{a_N^2}{\tilde{M}_N} \frac{1}{\tilde{M}_N} \sum_{i=1}^{\tilde{M}_N} f(\etpart{N,i}) \plim \beta \nu(f)
$$
and thus \eqref{eq:AN2} is satisfied.  \eqref{eq:AN3} is trivially satisfied.
\end{proof}

\begin{lem}
\label{lem:convRes}
Under the assumptions of Proposition \ref{thm:resid}, for any $f \in \consistfunc$,
$$
\frac{1}{\tilde{M}_N}\sum_{i=1}^{M_N} \left \lfloor \tilde{M}_N \sumweight{N}^{-1} \ewght{N,i}  \right\rfloor f(\epart{N,i})\plim
\nu \left( f \frac{\left\lfloor \ell \nu(1/\funcweight) \funcweight\right\rfloor}{\ell \; \nu(1/\funcweight) \funcweight} \right)\eqsp.
$$
\end{lem}

\begin{proof}
For any $K \geq 1$, denotes $\mathcal{B}_K= \bigcup_{j=0}^{\infty} [j-1/K,j+1/K]$. Because the weighted sample $\{ (\epart{N,i}, \funcweight(\epart{N,i}) \}_{1 \leq i \leq M_N}$
is consistent for $(\nu, \consistfunc)$ and  $\left \lfloor \tilde{M}_N \sumweight{N}^{-1} \ewght{N,i}  \right \rfloor \leq \tilde{M}_N  \sumweight{N}^{-1} \ewght{N,i}$,
we have for any $f \in \consistfunc$
\begin{align*}
& \tilde{M}_N^{-1}\sum_{i=1}^{M_N} \left \lfloor \tilde{M}_N \sumweight{N}^{-1} \ewght{N,i} \right \rfloor f(\epart{N,i})\1 \left\{ \ell \nu (1 / \funcweight) \ewght{N,i} \in (K,\infty)\cup \left([0,K]\cap
\mathcal{B}_K \right)\right\} \\
&\quad \leq \frac{1}{\sumweight{N}}\sum_{i=1}^{M_N} \ewght{N,i} f(\epart{N,i})\1 \left\{ \ell \nu(1/\funcweight) \funcweight(\epart{N,i}) \in (K,\infty)\cup \left([0,K]\cap \mathcal{B}_K
\right)\right\} \\
& \plim \int  f(\epart{}) \1 \left\{ \ell \nu(1/\funcweight) \funcweight(\epart{}) \in (K,\infty)\cup \left([0,K]\cap \mathcal{B}_K \right)\right\} \nu(d\epart{}) \eqsp.
\end{align*}
The \rhs\ of the previous display  can be made arbitrarily small by taking $K$ sufficiently because
$ \int f(\epart{}) \1 \left\{\ell \nu(1/\funcweight) \funcweight(\epart{}) \in \{\infty\} \cup \nset \right\} \nu(d\epart{})=0$.
For any $K \geq 1$, there exists  $\eta>0$ such that for any $a,b \in \rset$,
$$
\1 \{| a - \ell \nu(1/\funcweight)| \leq \eta\}\1 \left\{b \in [0,K] \setminus \mathcal{B}_K \right\} \left(\left \lfloor a b \right \rfloor -
\left \lfloor \ell \nu(1/\funcweight) b\right \rfloor\right)= 0
$$
Combining the previous equality with
$$
\tilde{M}_N \sumweight{N}^{-1}=\frac{\tilde{M}_N}{M_N} \left( \sumweight{N}^{-1} \sum_{i=1}^{M_N} \ewght{N,i}\frac{1}{\Phi(\epart{N,i})} \right) \plim \ell \nu(1/\funcweight)
$$
and
\begin{multline*}
 \tilde{M}_N^{-1}\sum_{i=1}^{M_N} \left\lfloor \ell \nu(1/\funcweight) \funcweight(\epart{N,i})\right\rfloor f(\epart{N,i})\1 \left\{ \ell \nu(1/\funcweight) \funcweight(\epart{}) \in [0,K] \setminus \mathcal{B}_K
 \right\}\\
\plim  \int \frac{\left\lfloor \ell \nu(1/\funcweight) \funcweight(\epart{})\right\rfloor}{\ell \; \nu(1/\funcweight) \funcweight(\epart{})}
f(\epart{})\1 \left\{ \ell \nu(1/\funcweight) \funcweight(\epart{}) \in [0,K] \setminus \mathcal{B}_K \right\} \nu(d\epart{})
\end{multline*}
yields
\begin{multline*}
\tilde{M}_N^{-1} \sum_{i=1}^{M_N} \left \lfloor \tilde{M}_N \sumweight{N}^{-1} \ewght{N,i}  \right \rfloor f(\epart{N,i})
\1 \left\{ \ell \nu(1/\funcweight) \funcweight(\epart{N,i}) \in [0,K] \setminus \mathcal{B}_K \right\} \\
\plim \int \frac{\left\lfloor \ell \nu(1/\funcweight) \funcweight(\epart{})\right\rfloor}{\ell \; \nu(1/\funcweight) \funcweight(\epart{})} f(\epart{})
\1 \left\{ \ell \nu(1/\funcweight) \funcweight(\epart{}) \in [0,K] \setminus \mathcal{B}_K\right\} \nu(d\epart{}) \eqsp.
\end{multline*}
The proof follows by letting $K\to \infty$.
\end{proof}

The condition $\nu\left\{\ell \nu(1/\Phi)\funcweight \in \nset \cup \{ \infty \} \right\}=0$ in Proposition \ref{thm:resid} and Lemma \ref{lem:convRes} is crucial.
Assume that $\{\epart{N,i}\}_{1 \leq i \leq N}$ is an \iid\ $\mu$-distributed sample where $\mu$ is the distribution on the set $\{1/2,2\}$ given by:
$\mu(\{1/2 \})=2/3$ and  $\mu(\{ 2 \})=1/3$. Let $\nu$ be the distribution on $\{1/2,2\}$ given by: $\nu(\{1/2\})=1/3$ and  $\nu(\{2\})=2/3$.
The weighted sample $\{ (\epart{N,i}, \epart{N,i}) \}_{1 \leq i \leq N}$ (\ie\ where we have set $\funcweight(\epart{})= \epart{}$) is a
consistent sample for $\nu$: for any function $f \in \mathbb{B}(\{1/2,2\}) \eqdef \{f: \ \{1/2,2\} \to \rset,  |f(1/2)|<\infty\quad \mbox{and} \quad |f(2)|<\infty \}$,
\begin{multline*}
  \frac{\sum_{i=1}^{M_N} \epart{N,i} f(\epart{N,i})}{\sum_{i=1}^{M_N} \epart{N,i}} \\
\plim \frac{(1/2) f(1/2) \mu(1/2) + 2 f(2) \mu(2)}{(1/2) \mu(1/2) + 2 \mu(2)} = (1/2) f(1/2) + 1/3 f(2)= \nu(f) \eqsp.
\end{multline*}
In this example, $\ell=1$ and obviously $\nu(1/\Phi)=1$. Moreover,
$$
\nu\left\{\funcweight \in \{\infty \}\cup \nset\right\}=\nu\left\{ \{1/2,2\} \cap \nset\right\}= \nu(\{2\})= 2/3\neq 0 \eqsp.
$$
We will show that the convergence in Lemma \ref{lem:convRes} fails. More precisely, setting $f(\epart{})=\epart{}$, we will show that
\begin{equation} \label{eq:counterEx}
\frac{1}{M_N}\sum_{i=1}^{M_N} \left \lfloor M_N \frac{\ewght{N,i}}{\sumweight{N}} \right\rfloor f(\epart{N,i}) =
\frac{1}{M_N}\sum_{i=1}^{M_N} \left \lfloor M_N \frac{\epart{N,i}}{\sum_{j=1}^{M_N} \epart{N,j}} \right\rfloor f(\epart{N,i}) \dlim   4/3 -(2Z)/3 \eqsp,
\end{equation}
where $Z$ is a Bernoulli variable with parameter $1/2$. This would imply that  $M_N^{-1} \sum_{i=1}^{M_N} \left \lfloor M_N \frac{\ewght{N,i}}{\sumweight{N}} \right\rfloor f(\epart{N,i})$
does not converge in probability to a constant. The LLN and CLT for \iid\ random variables imply that
\begin{align*}
&  M_N^{-1} \sumweight{N} = M_N^{-1} \sum_{i=1}^{M_N} \epart{N,i} \plim 1\\
& \begin{bmatrix} \1 \{M_N \sumweight{N}^{-1} <1 \}\\ \1\{M_N \sumweight{N}^{-1}\geq 1 \} \end{bmatrix} =
\begin{bmatrix}\1 \{(M_N)^{1/2}(M_N^{-1}\sumweight{N}-1)>0 \}\\ \1 \{(M_N)^{1/2}(M_N^{-1} \sumweight{N}-1)\leq 0 \}\end{bmatrix} \dlim \begin{bmatrix}Z \\ 1-Z\end{bmatrix}
\end{align*}
where $Z$ is a Bernoulli random variable with parameter $1/2$. Since $\ewght{N,i}=\funcweight(\epart{N,i})=\epart{N,i}  \in \{1/2,2\}$
and $f(\epart{})= \epart{}$,
\begin{align*}
&\1\left\{\frac{1}{2} <\frac{ M_N}{\sumweight{N}}<\frac{3}{2}\right\}\frac{1}{M_N}\sum_{i=1}^{M_N} \left \lfloor M_N \frac{\ewght{N,i}}{\sumweight{N}} \right\rfloor f(\epart{N,i}) \\
& \quad = \1\left\{\frac{1}{2} <\frac{ M_N}{\sumweight{N}}<\frac{3}{2}\right\}\frac{2}{M_N}\sum_{i=1}^{M_N} \left \lfloor  \frac{2M_N}{\sumweight{N}} \right\rfloor \1\{\epart{N,i}=2\}\\
&\quad = \1\left\{\frac{1}{2} <\frac{ M_N}{\sumweight{N}}<1\right\}\frac{2}{M_N} \sum_{i=1}^{M_N}  \1\{\epart{N,i}=2\}\\
& \quad \quad + \1\left\{1 \leq \frac{ M_N}{\sumweight{N}}<\frac{3}{2}\right\}\frac{4}{M_N} \sum_{i=1}^{M_N}  \1\{\epart{N,i}=2\}\\
& \quad \dlim (2Z)/3 +4(1-Z)/3=4/3-(2Z)/3\eqsp.
\end{align*}
The proof of (\ref{eq:counterEx}) is concluded by noting that $\frac{ M_N}{\sumweight{N}} \plim 1$.

%\bibliographystyle{ims}
%\bibliography{motherofallbibs}

\begin{thebibliography}{38}
\expandafter\ifx\csname natexlab\endcsname\relax\def\natexlab#1{#1}\fi
\expandafter\ifx\csname url\endcsname\relax
  \def\url#1{\texttt{#1}}\fi
\expandafter\ifx\csname urlprefix\endcsname\relax\def\urlprefix{URL }\fi

\bibitem[{Aldous and Eagleson(1978)}]{aldous:eagleson:1978}
\textsc{Aldous, D.~J.} and \textsc{Eagleson, G.~K.} (1978).
\newblock On mixing and stability of limit theorems.
\newblock \textit{Ann. Probab.} \textbf{6} 325--331.

\bibitem[{Berzuini and Gilks(2001)}]{berzuini:gilks:2001}
\textsc{Berzuini, C.} and \textsc{Gilks, W.~R.} (2001).
\newblock Resample-move filtering with cross-model jumps.
\newblock In \textit{Sequential {M}onte {C}arlo Methods in Practice}
  (A.~Doucet, N.~{De Freitas} and N.~Gordon, eds.). Springer.

\bibitem[{Capp\'e et~al.(2005)Capp\'e, Guillin, Marin and
  Robert}]{Cappe:Guillin:Marin:Robert:2003}
\textsc{Capp\'e, O.}, \textsc{Guillin, A.}, \textsc{Marin, J.} and
  \textsc{Robert, C.} (2005).
\newblock Population {M}onte {C}arlo.
\newblock \textit{J. Comput. Graph. Statist.} \textbf{13} 907--929.

\bibitem[{Chopin(2004)}]{chopin:2004}
\textsc{Chopin, N.} (2004).
\newblock Central limit theorem for sequential monte carlo methods and its
  application to bayesian inference.
\newblock \textit{Ann. Statist.} \textbf{32} 2385--2411.

\bibitem[{Crisan and Doucet(2002)}]{crisan:doucet:2002}
\textsc{Crisan, D.} and \textsc{Doucet, A.} (2002).
\newblock A survey of convergence results on particle filtering methods for
  practitioners.
\newblock \textit{IEEE Trans. Signal Process.} \textbf{50} 736--746.

\bibitem[{Crisan and Lyons(1997)}]{crisan:lyons:1997}
\textsc{Crisan, D.} and \textsc{Lyons, T.} (1997).
\newblock Nonlinear filtering and measure-valued processes.
\newblock \textit{Probab. Theory Related Fields} \textbf{109} 217--244.

\bibitem[{Cs{\"o}rg{\H{o}} and Fischler(1971)}]{csorgo:fischler:1971}
\textsc{Cs{\"o}rg{\H{o}}, M.} and \textsc{Fischler, R.} (1971).
\newblock On mixing and the central limit theorem.
\newblock \textit{T\^ohoku Math. J. (2)} \textbf{23} 139--145.

\bibitem[{{Del Moral}(1996)}]{delmoral:1996}
\textsc{{Del Moral}, P.} (1996).
\newblock Nonlinear filtering: interacting particle solution.
\newblock \textit{Markov Process. Related Fields} \textbf{2} 555--579.

\bibitem[{{Del Moral}(2004)}]{delmoral:2004}
\textsc{{Del Moral}, P.} (2004).
\newblock \textit{{F}eynman-Kac {F}ormulae. {G}enealogical and Interacting
  Particle Systems with Applications}.
\newblock Springer.

\bibitem[{Del~Moral and Guionnet(1999)}]{delmoral:guionnet:1999}
\textsc{Del~Moral, P.} and \textsc{Guionnet, A.} (1999).
\newblock Central limit theorem for nonlinear filtering and interacting
  particle systems.
\newblock \textit{Ann. Appl. Probab.} \textbf{9} 275--297.

\bibitem[{Del~Moral and Miclo(2000)}]{delmoral:miclo:2000}
\textsc{Del~Moral, P.} and \textsc{Miclo, L.} (2000).
\newblock Branching and interacting particle systems approximations of
  {F}eynman-{K}ac formulae with applications to non-linear filtering.
\newblock In \textit{S\'eminaire de Probabilit\'es, XXXIV}, vol. 1729 of
  \textit{Lecture Notes in Math.} Springer, Berlin, 1--145.

\bibitem[{Douc et~al.(2005)Douc, Guillin, Marin and
  Robert}]{douc:guillin:marin:robert:2004}
\textsc{Douc, R.}, \textsc{Guillin, A.}, \textsc{Marin, J.-M.} and
  \textsc{Robert, C.} (2005).
\newblock Convergence of adaptive sampling schemes.
\newblock Tech. Rep. 2005-6, CEREMADE.

\bibitem[{Doucet et~al.(2001)Doucet, {De Freitas} and
  Gordon}]{doucet:defreitas:gordon:2001}
\textsc{Doucet, A.}, \textsc{{De Freitas}, N.} and \textsc{Gordon, N.} (eds.)
  (2001).
\newblock \textit{Sequential {M}onte {C}arlo Methods in Practice}.
\newblock Springer, New York.

\bibitem[{Doucet et~al.(2000)Doucet, Godsill and
  Andrieu}]{doucet:godsill:andrieu:2000}
\textsc{Doucet, A.}, \textsc{Godsill, S.} and \textsc{Andrieu, C.} (2000).
\newblock On sequential {M}onte-{C}arlo sampling methods for {B}ayesian
  filtering.
\newblock \textit{Stat. Comput.} \textbf{10} 197--208.

\bibitem[{Dvoretzky(1972)}]{dvoretzky:1972}
\textsc{Dvoretzky, A.} (1972).
\newblock Asymptotic normality for sums of dependent random variables.
\newblock In \textit{Proceedings of the Sixth Berkeley Symposium on
  Mathematical Statistics and Probability (Univ. California, Berkeley, Calif.,
  1970/1971), Vol. II: Probability theory}. Univ. California Press, Berkeley,
  Calif.

\bibitem[{Gelman and Rubin(1992)}]{gelman:rubin:1992}
\textsc{Gelman, A.} and \textsc{Rubin, D.~B.} (1992).
\newblock Inference from iterative simulation using multiple sequences.
\newblock \textit{Statist. Sci.} \textbf{7} 473--483.

\bibitem[{Gilks and Berzuini(2001)}]{gilks:berzuini:2001b}
\textsc{Gilks, W.~R.} and \textsc{Berzuini, C.} (2001).
\newblock Following a moving target---{M}onte {C}arlo inference for dynamic
  {B}ayesian models.
\newblock \textit{J. Roy. Statist. Soc. Ser. B} \textbf{63} 127--146.

\bibitem[{Hall and Heyde(1980)}]{hall:heyde:1980}
\textsc{Hall, P.} and \textsc{Heyde, C.} (1980).
\newblock \textit{Martingale Limit Theory and its Application}.
\newblock Academic Press, New York, London.

\bibitem[{Hall and Heyde(1981)}]{hall:heyde:1981}
\textsc{Hall, P.} and \textsc{Heyde, C.~C.} (1981).
\newblock Rates of convergence in the martingale central limit theorem.
\newblock \textit{Ann. Probab.} \textbf{9} 395--404.

\bibitem[{Handschin and Mayne(1969)}]{handschin:mayne:1969}
\textsc{Handschin, J.} and \textsc{Mayne, D.} (1969).
\newblock {M}onte {C}arlo techniques to estimate the conditionnal expectation
  in multi-stage non-linear filtering.
\newblock In \textit{Int. J. Control}, vol.~9.

\bibitem[{Kong et~al.(1994)Kong, Liu and Wong}]{kong:liu:wong:1994}
\textsc{Kong, A.}, \textsc{Liu, J.~S.} and \textsc{Wong, W.} (1994).
\newblock Sequential imputation and {B}ayesian missing data problems.
\newblock \textit{J. Am. Statist. Assoc.} \textbf{89} 590--599.

\bibitem[{K\"{u}nsch(2001)}]{kuensch:2001}
\textsc{K\"{u}nsch, H.~R.} (2001).
\newblock State space and hidden markov models.
\newblock In \textit{Complex Stochastic Systems} (O.~E. Barndorff-Nielsen,
  D.~R. Cox and C.~Klueppelberg, eds.). CRC Publisher, Boca raton, 109--173.

\bibitem[{K\"{u}nsch(2003)}]{kuensch:2003}
\textsc{K\"{u}nsch, H.~R.} (2003).
\newblock Recursive {M}onte-{C}arlo filters: algorithms and theoretical
  analysis.
\newblock Preprint ETHZ, seminar f\"{u}r statistics.

\bibitem[{Landau and Binder(2000)}]{landau:binder:2000}
\textsc{Landau, D.~P.} and \textsc{Binder, K.} (2000).
\newblock \textit{A guide to {M}onte {C}arlo simulations in statistical
  physics}.
\newblock Cambridge University Press, Cambridge.

\bibitem[{Liu(2001)}]{liu:2001}
\textsc{Liu, J.} (2001).
\newblock \textit{{M}onte {C}arlo Strategies in Scientific Computing}.
\newblock Springer, New York.

\bibitem[{Liu and Chen(1995)}]{liu:chen:1995}
\textsc{Liu, J.} and \textsc{Chen, R.} (1995).
\newblock Blind deconvolution via sequential imputations.
\newblock \textit{J. Roy. Statist. Soc. Ser. B} \textbf{430} 567--576.

\bibitem[{Liu and Chen(1998)}]{liu:chen:1998}
\textsc{Liu, J.} and \textsc{Chen, R.} (1998).
\newblock Sequential {M}onte-{C}arlo methods for dynamic systems.
\newblock \textit{J. Roy. Statist. Soc. Ser. B} \textbf{93} 1032--1044.

\bibitem[{Liu et~al.(2001)Liu, Chen and Logvinenko}]{liu:chen:logvinenko:2001}
\textsc{Liu, J.}, \textsc{Chen, R.} and \textsc{Logvinenko, T.} (2001).
\newblock A theoretical framework for sequential importance sampling and
  resampling.
\newblock In \textit{Sequential {M}onte {C}arlo Methods in Practice}
  (A.~Doucet, N.~{De Freitas} and N.~Gordon, eds.). Springer.

\bibitem[{Liu(1996)}]{liu:1996}
\textsc{Liu, J.~S.} (1996).
\newblock Metropolized independent sampling with comparisons to rejection
  sampling and importance sampling.
\newblock \textit{Stat. Comput.} \textbf{6} 113--119.

\bibitem[{McLeish(1974)}]{mcleish:1974}
\textsc{McLeish, D.~L.} (1974).
\newblock Dependent central limit theorems and invariance principles.
\newblock \textit{Ann. Probab.} \textbf{2} 620--628.

\bibitem[{Petrov(1995)}]{petrov:1995}
\textsc{Petrov, V.~V.} (1995).
\newblock \textit{Limit Theorems of Probability Theory}.
\newblock Oxford University Press.

\bibitem[{Pitt and Shephard(1999)}]{pitt:shephard:1999}
\textsc{Pitt, M.~K.} and \textsc{Shephard, N.} (1999).
\newblock Filtering via simulation: Auxiliary particle filters.
\newblock \textit{J. Am. Statist. Assoc.} \textbf{94} 590--599.

\bibitem[{Ristic et~al.(2004)Ristic, Arulampalam and
  Gordon}]{ristic:arulampalam:gordon:2004}
\textsc{Ristic, B.}, \textsc{Arulampalam, M.} and \textsc{Gordon, A.} (2004).
\newblock \textit{Beyond Kalman Filters: Particle Filters for Target Tracking}.
\newblock Artech House.

\bibitem[{Rubin(1987)}]{rubin:1987}
\textsc{Rubin, D.~B.} (1987).
\newblock A noniterative sampling/importance resampling alternative to the data
  augmentation algorithm for creating a few imputations when the fraction of
  missing information is modest: the {SIR} algorithm (discussion of {T}anner
  and {W}ong).
\newblock \textit{J. Am. Statist. Assoc.} \textbf{82} 543--546.

\bibitem[{Shephard and Pitt(1997)}]{shephard:pitt:1997}
\textsc{Shephard, N.} and \textsc{Pitt, M.} (1997).
\newblock Likelihood analysis of non-{G}aussian measurement time series.
\newblock \textit{Biometrika} \textbf{84} 653--667.
\newblock Erratum in 91:249--250, 2004.

\bibitem[{Shiryaev(1996)}]{shiryaev:1996}
\textsc{Shiryaev, A.~N.} (1996).
\newblock \textit{Probability}.
\newblock 2nd ed. Springer.

\bibitem[{Tanizaki(1999)}]{tanizaki:1999}
\textsc{Tanizaki, H.} (1999).
\newblock On the nonlinear and nonnormal filter using rejection sampling.
\newblock \textit{IEEE Trans. Automat. Control} \textbf{44} 314--319.

\bibitem[{Tanizaki(2001)}]{tanizaki:2001}
\textsc{Tanizaki, H.} (2001).
\newblock Nonlinear and non-{G}aussian state space modeling using sampling
  techniques.
\newblock \textit{Ann. Inst. Statist. Math.} \textbf{53} 63--81.

\end{thebibliography}

\end{document}